\documentclass[12pt]{article}
\usepackage[utf8]{inputenc}
\usepackage{amssymb}
\usepackage{amsmath}
\usepackage{amsthm}
\usepackage{hyperref}
\usepackage{a4wide}
\usepackage{breqn}
\usepackage[table,xcdraw]{xcolor}
\usepackage{xcolor}
\usepackage{mathtools}
\usepackage{bm}
\usepackage{enumitem}
\usepackage{mathrsfs}
\usepackage{graphicx}
\usepackage{caption}
\usepackage{url}
\usepackage{euscript}
\usepackage{placeins}
\usepackage{mathdots}
\usepackage[english]{babel}
\usepackage{blkarray}
\usepackage{array}
\usepackage{subcaption}
\usepackage{booktabs}
\usepackage[numbers]{natbib}

\usepackage{graphics}
\usepackage{natbib}
\newcommand{\indsize}{\scriptsize}
\newcommand{\colind}[2]{\displaystyle\smash{\mathop{#1}^{\raisebox{.5\normalbaselineskip}{\indsize #2}}}}
\newcommand{\rowind}[2]{\displaystyle\smash{\mathop{#1}_{\raisebox{-.5\normalbaselineskip}{\indsize #2}}}}
\numberwithin{equation}{section}

\newtheorem{theorem}{\bf Theorem}[section]
\newtheorem{definition}{Definition}[section]
\newtheorem{corollary}{Corollary}[section]

\newtheorem{problem}{Problem}[section]
\newtheorem{lemma}{Lemma}[section]
\newtheorem{remark}{Remark}[section]
\theoremstyle{remark}
\newtheorem{exam}{\bf Example}

\def \vec{\mathrm v\mathrm e \mathrm c}

\def \R{{\mathbb R}}

\def \C{{\mathbb{C}}}
\def \x{\widetilde x}
\def \y{\widetilde y}

\def\bmatrix#1{\left[\begin{matrix}
		#1
	\end{matrix}\right]}
\def \diag{\mathrm{diag}}

\def \E{{\mathcal E}}

\def \CR{{\mathcal C}}
\def \A{{\mathcal A}}
\def \ST{{\mathcal S \mathcal T}}
\def \T{\mathcal T}
\def \D{{\Delta}}
\def \a{{\bf a}}

\def \m{{\bf m}}

\def \g{{\bf g}}
\def \t{{\bf t}}
\def \r{{\bf r}}
\def \bdot{{\, \odot\, }}
\newcommand{\vertiii}[1]{{\left\vert\kern-0.25ex\left\vert\kern-0.25ex\left\vert #1 
    \right\vert\kern-0.25ex\right\vert\kern-0.25ex\right\vert}}

\def \rank{\mathrm{rank}}

\def\bmatrix#1{\left[ \setlength{\arrayrulewidth}{20pt}
\begin{matrix} #1 \end{matrix} \right]}

\def \R{{\mathbb R}}

\newcommand{\Proof}{\textit{Proof. \,}}
\date{}

\title{
Structured Backward Errors for Special Classes of Saddle Point Problems with Applications}
\author{Sk. Safique Ahmad\footnotemark[2]  \and Pinki Khatun\footnotemark[2] }

\begin{document}
	
	\maketitle
	\begin{abstract}
  {In the realm of numerical analysis, the study of structured backward errors (\textit{BEs}) in saddle point problems (\textit{SPPs}) has shown promising potential for development.}
   However, these investigations overlook the inherent sparsity pattern of the coefficient matrix of the \textit{SPP}. Moreover, the existing techniques are not applicable when the block matrices have  {\textit{circulant}}, \textit{Toeplitz}, or \textit{symmetric}-\textit{Toeplitz} structures and do not even provide structure preserving minimal perturbation matrices for which the \textit{BE} is attained. To overcome these limitations, we investigate the structured \textit{BEs} of
   \textit{SPPs} when the perturbation matrices exploit the sparsity pattern as well as  {\textit{circulant}}, \textit{Toeplitz}, and \textit{symmetric}-\textit{Toeplitz} structures. Furthermore, we construct minimal perturbation matrices that preserve the sparsity pattern and the aforementioned structures.  {Applications of the developed frameworks are utilized to compute \textit{BEs} for the weighted regularized least squares problem.} 
Finally, numerical experiments are performed to validate our findings, showcasing the utility of the obtained structured \textit{BEs} in assessing the strong backward stability of numerical algorithms.
	\end{abstract}
 
	
	\footnotetext[2]{
		Department of Mathematics, Indian Institute of Technology Indore,  Khandwa Road, Simrol, Indore 453552, Madhya Pradesh, India, \texttt{safique@iiti.ac.in}, \texttt{pinki996.pk@gmail.com}.}
	\section{Introduction}
 Saddle point problems (\textit{SPPs}) have garnered significant attention in recent times due to their pervasive occurrence across diverse domains of scientific computation and engineering applications \cite{Benzi2005}. These domains encompass constrained optimization \cite{OPTM1992, OPTM2014}, computational fluid dynamics \cite{Elman2005}, constrained and weighted least square estimation \cite{LSproblem2002}, and linear elasticity \cite{SymbolCirculant}. The standard \textit{SPP} is a  two-by-two block linear system of the following  form:
\begin{align}\label{eq11}
    \mathcal{A}{u}\triangleq \bmatrix{A&B^T\\B &C}\bmatrix{x\\y}=\bmatrix{f\\g}\triangleq d,
\end{align}
where $A\in \C^{n\times n}, B\in \C^{m\times n}, C\in \C^{m\times m},f\in \C^n,$ and $g\in \C^m.$ Here, $B^T$ stands for transpose of $B$ and {$\C^{m\times n}$ denotes the set of all $m\times  n$ complex matrices}. In general, the block matrices $A,B$ and $C$ are sparse \cite{Elman2005}.
Further, the block matrices in \eqref{eq11} can have \textit{symmetric}, \textit{Toeplitz}, \textit{symmetric}-\textit{Toeplitz}, or  {\textit{circulant}} structures.
For instance, consider the following weighted  regularized least squares (\textit{WRLS}) problem  \cite{circulantbased2021}: 
\begin{equation}\label{ls:eq12}
    \min_{z\in \R^{n}}\| {M}z-b\|_2,
 \end{equation} where $ {M}=\bmatrix{W^{\frac{1}{2}}K^{T}\\\sqrt{\lambda}I_m}$ and $b=\bmatrix{W^{\frac{1}{2}}f\\ {\bf 0}}$ with $W\in \R^{n\times n}$ is a \textit{symmetric} positive definite weighting matrix, $ {K\in \R^{m\times n}}$ is a \textit{Toeplitz} or \textit{symmetric}-\textit{Toeplitz}  matrix $ {(n=m)}$, $I_m$ is the identity matrix of size $m$ and $\lambda >0$ is the regularization parameter. {Here, $\R^{m\times n}$ denotes the set of all $m\times  n$ real matrices}. Problem \eqref{ls:eq12} plays a vital role in image reconstruction \cite{imagereconstruction} and image restoration with colored noise \cite{image1989}. While solving \eqref{ls:eq12}, we encounter \textit{SPPs} of the form \eqref{eq11}.
 Also, \textit{SPP} involving  {\textit{circulant}} or \textit{Toeplitz} block matrices often arise during the discretization of elasticity problems using finite difference scheme \cite{SymbolCirculant}.

 In recent times, a number of numerical algorithms have been developed to find the efficient solution of the \textit{SPP} \eqref{eq11} with  {\textit{circulant}}, \textit{Toeplitz}, or \textit{symmetric}-\textit{Toeplitz} block matrices; see \cite{weightedLS2023, SymbolCirculant, circulantbased2021, SYMTOEP2017}. 
 This prompts a natural inquiry: can an approximate solution obtained from a numerical algorithm serve as the exact solution to a nearly perturbed problem?  {The concept of \textit{BE} is used to determine how far a computed solution stands from the original problem.  By multiplying the condition number with the \textit{BE}, an upper bound on the forward error can be established \cite{higham2002}.} 
Moreover, \textit{BE} serves as an effective and reliable stopping criterion for numerical algorithms solving large system of linear equations \cite{stopping2009}.

 Assume that $ { {\widetilde{u}}}=[\widetilde{x}^{T},\widetilde{y}^{T}]^{T}$ is  an approximate solution of the system \eqref{eq11}. Then, the normwise unstructured backward error $\eta( { {\widetilde{u}}})$ is defined as: 
\begin{equation*}
\eta^{\mathcal{S}_i}(\widetilde{x},\widetilde{y})=\displaystyle{ \min_{\substack{ \D\mathcal{A}\in\, \C^{(n+m)\times (n+m)}, \, \D d\in \C^{n+m},\\ (\A+\D \A)\widetilde{u}=d+\D d}}} \left\|\bmatrix{\frac{\|\Delta \A\|_F}{\|\A\|_F}, &  {\frac{\|\Delta {d}\|_2}{\|{d}\|_2}}}\right\|_2,
     \end{equation*}
which has a compact formula, given by \citet{Rigal1967} as follows: 
\begin{equation}\label{eq14}
    \eta( { {\widetilde{u}}})=\frac{\|{d}-\A  { {\widetilde{u}}}\|_2}{\sqrt{\|\A\|^2_F\| { {\widetilde{u}}}\|_2^2+\|{d}\|_2^2}}.
\end{equation}
Here, $\|\cdot\|_2$ and $\|\cdot\|_F$ stand for the Euclidean and Frobenius norms, respectively. A small value of $\eta( {\widetilde{u}})$ indicates that the approximate solution $ {\widetilde{u}}$ is the exact solution of a slightly perturbed system $(\A +\Delta \A) {\widetilde{u}}={d}+\Delta {d}$ with relatively small $\|\D \A\|_F$ and $\|\D {d}\|_2$.  {A numerical algorithm to solve a problem is backward stable if the computed solution to the problem is always the exact solution of a nearby problem \cite{higham2002}}. 
Notwithstanding, the perturbed coefficient matrix $\A +\Delta \A$ does not necessarily follow the structure of \eqref{eq11}. An algorithm is said to be strongly backward stable \cite{strongweak, strongstab}  {if the computed solution is always the exact solution of a nearby structure-preserving problem.} 
Consequently, it becomes intriguing to examine whether preserving the special structure of the block matrices of $\A$ in the perturbation matrix $\Delta \A$ leads to small  $\|\Delta \A\|_F$ or not. That is, the numerical algorithms for solving \textit{SPPs} are strongly backward stable or not. So far, significant research has been done in this direction, and structured \textit{BE} has been extensively considered in \cite{Sun1999, be2007wei, be2012LAA, be2017ma, BE2020BING, be2022lma}, where the perturbation matrix  $\Delta \A$ preserve the original structure of $\A$ and perturbation on block matrices $A$ or $C$ preserve the \textit{symmetric} structures.  However, it is noteworthy to mention a few drawbacks of the aforementioned studies: 
\begin{itemize}
    \item The coefficient matrix $\A$ of the \textit{SPP} \eqref{eq11} is generally sparse, and the existing studies do not consider and preserve the sparsity pattern of the coefficient matrix $\A.$ 
    \item Existing techniques are not applicable when the block matrices $A$, $B$ and $C$ in \textit{SPP} \eqref{eq11} have  {\textit{circulant}}, \textit{Toeplitz}, or \textit{symmetric}-\textit{Toeplitz} structures.
    \item Moreover,  the research available in the literature for structured \textit{BE} analysis for \eqref{eq11} does not provide the explicit formulae for the minimal perturbation matrices for which structured \textit{BE} is attained and preserves the inherent matrix structure.
\end{itemize}
 
By preserving the sparsity pattern of the original matrices, structured \textit{BEs} have been studied in the literature;
 see, for example,  \cite{prince2020, prince2021, PEVP2012}. The block matrices in \eqref{eq11} are often sparse in many applications, making it essential to maintain their sparsity pattern in the perturbation matrices. This paper addresses the aforementioned challenges by investigating structured \textit{BEs} for the \textit{SPP} \eqref{eq11} by preserving both the inherent block structure and sparsity in the perturbation matrices under three scenarios. First, we consider the case where \(n = m\) and the block matrices \(A\), \(B\), and \(C\) are \textit{circulant}. Second, we consider \(A\), \(B\), and \(C\) are \textit{Toeplitz} matrices. Third, we analyze the case when \(n = m\), \(B\) is \textit{symmetric}-\textit{Toeplitz}, and \(A, C \in \mathbb{C}^{n \times n}\). The following are the main contributions of this paper: 
\begin{itemize}
    \item We investigate the structured \textit{BEs} for the \textit{SPP} \eqref{eq11} when the block matrices $A, B$ and $C$  {possesses}  {\textit{circulant}}, \textit{Toeplitz}, and \textit{symmetric}-\textit{Toeplitz} structure with or without preserving the sparsity pattern.
    \item We develop frameworks that give the minimal perturbation matrices that retain the  {\textit{circulant}}, \textit{Toeplitz}, or \textit{symmetric}-\textit{Toeplitz} structures, as well as the sparsity pattern of the original matrices. 
    \item  {We provide an application of our obtained results in finding the structured \textit{BEs} for the \textit{WRLS} problem when the coefficient matrix exhibits  \textit{Toeplitz} or \textit{symmetric}-\textit{Toeplitz} structure}. 
    \item Lastly, numerical experiments are performed to test the backward stability and strong backward stability of numerical algorithms to solve \textit{SPP} \eqref{eq11}.
\end{itemize}
This paper is organized as follows. In Section \ref{SEC2}, we discuss some notations, preliminary definitions, and results. In Sections \ref{SEC3}-\ref{SEC5}, structured \textit{BEs} for  {\textit{circulant}}, \textit{Toeplitz}, and \textit{symmetric}-\textit{Toeplitz} matrices are derived, respectively. Moreover, in Section \ref{SEC5}, we discuss the unstructured \textit{BE} for the \textit{SPP} \eqref{eq11} by only preserving the sparsity pattern. In Section \ref{SEC6}, we provide an application of our developed theories. Section \ref{SEC:Numerical} includes a few numerical experiments. Finally, Section \eqref{SEC:Conclusion} provides some concluding remarks.
 
 \section{Notation and Preliminaries}\label{SEC2}
 Throughout the paper, we adopt the following notations. 
 {We use} $\CR_n,  {\T_{m\times n}}$ and $\ST_n$ to denote the collections of all complex $n\times n$   {\textit{circulant}},  {$m\times n$ \textit{Toeplitz}} and  {$n\times n$ } \textit{symmetric}-\textit{Toeplitz} matrices, respectively. For any matrix $A\in \C^{m\times n},$  $A^{\dagger}$ and $A^H$ stand for  the Moore-Penrose inverse and conjugate transpose of $A,$ respectively.  For  $A=[a_{ij}]\in \C^{m\times n},$ $\vec(A):=[{\bm{a}}_1^T,{\bm{a}}_2^T,\ldots,{\bm{a}}_n^T]^T\in \C^{mn},$ where ${\bm{a}}_i$ is the $i$-th column of  $A.$  {The notation} ${\bf 1}_{m\times n}\in \R^{m\times n}$ denotes the matrix with all entries are equal to $1.$ The  {sparsity pattern} of a matrix  $A\in \C^{m\times n}$ is defined as $\Theta_A:={\tt sgn}(A)=[{\tt sgn}(a_{ij})],$ where
 \begin{align*}
     {\tt sgn}(a_{ij})=\left\{ \begin{array}{lcl}
     1, &a_{ij}\neq 0,\\
     0,  &a_{ij}=0.
      \end{array}\right.
 \end{align*}
 The Hadamard product of $A,B\in \C^{m\times n}$ is defined as $A\odot B=[a_{ij}b_{ij}]\in \C^{m\times n}.$ For any vector $x\in \C^m,$ $\mathfrak{D}_{x}$ is the diagonal matrix defined as $\mathfrak{D}_{{x}}:=\diag({x})\in\C^{m\times m}.$ We denote ${\bf 0}_{m\times n}$ as the zero matrix of size $m\times n$ (for simplicity, we use ${\bf 0}$ when  matrix size is clear).  Let ${ w}:=[w_1, w_2, w_3, w_4, w_5]^T$, where $w_i$  are nonnegative real numbers for $i=1,2,\ldots,5,$ with the convention that $w_i^{-1}=0,$ whenever $w_i=0.$ For any ${ w},$ we define 
$$\vertiii{\bmatrix{\A& {d}}}_{{ w},F}=\Big\|\bmatrix{w_1\| A\|_F,& w_2\| B\|_F,& w_3\| C\|_F,& w_4\|f\|_2,&w_5  \|  g\|_2}\Big\|_2.$$
Note that $w_i=0$ implies that the corresponding block matrix has no perturbation. Next, we recall the definitions of  {\textit{circulant}}, \textit{Toeplitz}, and \textit{symmetric}-\textit{Toeplitz} matrices.
 
 \begin{definition}\label{def21}
     A matrix $C\in \mathbb{C}^{n\times n}$ is called a  {\textit{circulant}} matrix if for any vector $c=[c_1,c_2,\ldots,c_n]^T\in \mathbb{C}^n,$ it has the following form:
     \begin{equation}\label{eq21}
        {\tt Cr}(c):= C=\bmatrix{c_1 & c_n & c_{n-1}  &\cdots &c_2\\
         c_2 & c_1 & c_{n} & \cdots &c_3\\
         \vdots & \ddots&\ddots &\ddots& \vdots\\
         c_{n-1}&c_{n-2}&\ddots&\ddots&c_n\\
          c_n & c_{n-1} & \cdots & c_2& c_1}.
     \end{equation}
 \end{definition}
  We denote the generator vector for the  {\textit{circulant}} matrix $C$ as $$\vec_{\CR}(C):=[c_1,c_2,\ldots,c_n]^T\in \mathbb{C}^n.$$
 \begin{definition}\label{def22}
     A matrix $T=[t_{ij}]\in \mathbb{C}^{ {m}\times n}$ is called a \textit{Toeplitz} matrix if for any vector $$\vec_{\mathcal{T}}(T):=[t_{- {m}+1},t_{- {m}+2},\ldots,t_{-1},t_0,t_1,\ldots,t_{n-1}]^T\in \mathbb{C}^{ {n+m-1}},$$ we have $t_{ij}=t_{j-i},$ for all  {$1\leq i\leq m$ and $1\leq j\leq n$}.
 \end{definition}
 We denote $\vec_{\mathcal{T}}(T)$ as the generator vector for the \textit{Toeplitz} matrix $T.$ Also, for any vector $t\in \mathbb{C}^{ m+n-1}$  corresponding generated \textit{Toeplitz} matrix  is denoted by ${\mathcal{T}}(t).$
\begin{remark}
   The Toeplitz matrix $T$    is known as a \textit{symmetric}-\textit{Toeplitz} matrix when $ {n=m}$ and  $t_{-m+1}=t_{n-1},\ldots, t_{-1}=t_1.$ In this context, we employ the notation $$\vec_{\ST}(T):=[t_0,\ldots,t_{n-1}]^T\in \C^n$$ to denote its generator vector. Also, for any vector $t\in \C^n,$ the corresponding \textit{symmetric}-\textit{Toeplitz} matrix is symbolized as ${\ST}(t)$.
\end{remark}

Throughout the paper, we assume that the coefficient matrix $\A$ in \eqref{eq11}  is nonsingular. If the block matrices $A, B,$  and $C$ have  {\textit{circulant}} (or \textit{Toeplitz}) structure, we identify \eqref{eq11} as  {\textit{circulant}} (or \textit{Toeplitz}) structured \textit{SPP}.  Moreover, we call \eqref{eq11} as \textit{symmetric}-\textit{Toeplitz} structured \textit{SPP} when $B\in \ST_n$ and $A,C \in\C^{n\times n}.$

Next, we define normwise structured \textit{BE} for the \textit{SPP} \eqref{eq11}. Before this, we denote $$\D \A:=\bmatrix{\D A&\D B^T\\ \D B& \D C}  ~\text{and}~ \D {d}:=\bmatrix{\D f\\ \D g}.$$
\begin{definition}\label{def24}
     Let $ {\widetilde{u}}=[\widetilde{x}^T,\widetilde{y}^T]^T$ be an approximate solution of the \textit{SPP} \eqref{eq11}. Then,  the  normwise structured \textit{BEs} are defined as follows:
 { \begin{equation*}
\eta^{\mathcal{S}_i}(\widetilde{x},\widetilde{y})=\displaystyle{\min_{\substack{\D \mathcal{A}\in\, \mathcal{S}_i, \, \D d\in \C^{n+m},\\ (\A+\D \A)\widetilde{u}=d+\D d}}} \vertiii{\bmatrix{\D \A& \D{d}}}_{{ w},F},
     \quad \mbox{for}\, \, i=1,2,3,
     \end{equation*}  
     where
    $ \mathcal{S}_{1}=\Big\{\D \A \,{\big |} \, \D A,\,\D B,\, \D C\in \CR_{n}\Big\}, ~~
      \mathcal{S}_2=\Big\{\D \A \,{\big |}\, \D A\in \mathcal{T}
_{n\times n},\,\D B\in \mathcal{T}
_{m\times n},\, \D C\in \mathcal{T}
_{m\times m}\Big\}~
     \text{and}~ \mathcal{S}_3=\Big\{ \D\A~ {\big|}\, \D B\in \mathcal{ST}_n,\,\D A, \D C\in \C^{n\times n}\Big\}.$
 }
\end{definition}
 \begin{problem}\label{pr1}
      Find out the minimal  perturbation matrices $$ {\widehat{\D \A}:=\bmatrix{\widehat{\D A}&\widehat{\D B}^T\\ \widehat{\D B}& \widehat{\D C}}\in \mathcal{S}_{i}   
     ~ ~ \text{and}}~  {{\widehat{\D d}}:=\bmatrix{\widehat{\D f}\\ \widehat{\D g}}}$$
     such that
     $\eta^{\mathcal{S}_i}(\widetilde{x},\widetilde{y})= \vertiii{\bmatrix{\widehat{\D\A}& {\widehat{\D d}}}}_{{ w},F},~\text{for}~ i=1,2,3.$
 \end{problem}
\begin{remark}
Our main focus is studying perturbations with the same sparsity pattern as the original matrices. To achieve this,  we replace the perturbation matrices $\D A,\D B$ and $\D C$ by $\D A\bdot \Theta_A, \D B\bdot \Theta_B$ and $\D C\bdot \Theta_C,$ respectively. In this context, we denote the structured \textit{BEs} by $\eta_{\tt sps}^{\mathcal{S}_i}(\x,\y),$ $i=1,2,3.$ Further, the minimal perturbation matrices are denoted by $ \widehat{\D A}_{\tt sps},$ $\widehat{\D B}_{\tt sps},$ $\widehat{\D C}_{\tt sps},$ $\widehat{\D f}_{\tt sps},$ and  $\widehat{\D g}_{\tt sps}.$ Note that, if $M\in \CR_n\, (or\, \mathcal{T}_{ m\times n}, \, or \, \mathcal{ST}_n),$ then $\Theta_M\in \CR_n~(or\, \mathcal{T}_{ m\times n}, \, or \, \mathcal{ST}_n). $
\end{remark}
 

 \section{ Structured \textit{BEs} 
 for  {\textit{circulant}} structured \textit{SPPs}}\label{SEC3}
 In this section, we  {consider $n=m$ and}  derive explicit formulae for the structured \textit{BEs} $ \eta^{\mathcal{S}_1}_{\tt sps}(\x,\y)$ and $ \eta^{\mathcal{S}_1}(\x,\y),$ by preserving the  {\textit{circulant}} structure to the perturbation matrices. Moreover, we provide minimal perturbation matrices to the  {Problem \eqref{pr1}}. In order to obtain structured \textit{BEs} formulae, we derive the following lemma.
 \begin{lemma}\label{Lm31}
     Let $A,B,M\in\CR_n$ with  generator vectors $\vec_{\CR}(A)=[a_1,\ldots,a_n]^T\in \C^n,\vec_{\CR}(B)=[b_1,\ldots,b_n]^T\in \C^n,$ and $\vec_{\CR}(M)=[m_1,\ldots,m_n]^T\in \C^n,$ respectively. Suppose $x=[x_1,\ldots,x_n]^T\in \C^n$ and   $y=[y_1,\ldots,y_n]^T\in \C^n.$ 
     Then 
     { \begin{align*}
         &(A\,\odot\,\Theta_{M})x={\tt Cr}\left(x\right) \mathfrak{D}_{c(M)}\vec_{\CR}(A\,\odot\,\Theta_{M})~ \text{and}\\
         & 
         (B\,\odot\,\Theta_{M})^Ty=\mathcal{H}_y\mathfrak{D}_{c(M)}\vec_{\CR}(B\,\odot\,\Theta_{M}),
     \end{align*}}
      where  $c(M): =\vec_{\CR}(\Theta_M)$ and
 $\mathcal{H}_y\in\C^{n\times n}$ has the following form:
  \begin{equation}\label{eq210}
           \mathcal{H}_y=\bmatrix{y_1 & y_2 & \cdots & y_{n-1}& y_n\\ 
            y_2 & y_3 &\cdots & y_n & y_1\\
              \vdots &  \iddots  & \iddots& \iddots&\vdots\\
                 y_{n-1} & \iddots  &\iddots& y_{n-3}  & y_{n-2}\\
                  y_n & y_1  & \cdots& y_{n-2} & y_{n-1}}.
                                            \end{equation}      
 \end{lemma}
 \proof
   Since  $ij$-th entry of $A\bdot \Theta_M$ is $(A\bdot \Theta_M)_{ij}=a_{ij}\,{\tt sgn}(m_{ij}),$  we get $A\,\odot\,\Theta_{M}\in \CR_n,$ and $$\vec_{\CR}(A\,\odot\,\Theta_{M})=\bmatrix{a_1{\tt sgn}(m_1)\\ \vdots\\ a_n{\tt sgn}(m_n)}.$$ Now, expanding $(A\,\odot\,\Theta_{M})x,$ we get the following: 
    {\begin{align*}
    (A\,\odot\,\Theta_{M})x=  \bmatrix{ a_1\, {\tt sgn}(m_1)x_1+ a_n\, {\tt sgn}(m_n)x_2+\cdots+ a_2\, {\tt sgn}(m_2)x_n\\
       a_2\, {\tt sgn}(m_2)x_1+ a_1\, {\tt sgn}(m_1)x_2+\cdots+ a_3\, {\tt sgn}(m_3)x_n\\
       \quad \quad\vdots \hspace{2cm}\quad \vdots \hspace{2cm}\, \,\vdots \hspace{2cm} \vdots \\
       a_n\, {\tt sgn}(m_n)x_1+ a_{n-1}\, {\tt sgn}(m_{n-1})x_2+\cdots+ a_1\, {\tt sgn}(m_1)x_n}.
   \end{align*}}
  Since $({\tt sgn}(m_i))^2={\tt sgn}(m_i),$ rearrangement of the above gives 
  {\footnotesize \begin{align*}
     (A\,\odot\,\Theta_{M})x=\bmatrix{x_1\, {\tt sgn}(m_1)& x_n\, {\tt sgn}(m_2)& x_{n-1}\, {\tt sgn}(m_{3})& \cdots & x_2\, {\tt sgn}(m_n)\\  x_2\, {\tt sgn}(m_1)& x_1\, {\tt sgn}(m_2)& x_n\, {\tt sgn}(m_3)&\cdots& x_3\, {\tt sgn}(m_n)\\
    \vdots&\ddots& \ddots &\ddots& \vdots\\
     x_{n-1}\, {\tt sgn}(m_1)& x_{n-2}\, {\tt sgn}(m_2)&\cdots&\ddots& x_n\, {\tt sgn}(m_n)\\
      x_n\, {\tt sgn}(m_1)&  x_{n-1}\, {\tt sgn}(m_2) &\cdots& x_2\, {\tt sgn}(m_{n-1})&  x_1\, {\tt sgn}(m_n)}  \bmatrix{a_1{\tt sgn}(m_1)\\a_2{\tt sgn}(m_2)\\ \vdots \\
      \vdots \\a_n{\tt sgn}(m_n)}.
  \end{align*}}
Hence, the above can be expressed as
$$ {(A\,\odot\,\Theta_{M})x=} ~{\tt Cr}(x)\mathfrak{D}_{c(M)}\vec_{\CR}(A\,\odot\,\Theta_{M}).$$

Similarly, expanding $(B\,\odot\,\Theta_{M})^Ty,$ we can obtain 
   $$  {(B\,\odot\,\Theta_{M})^Ty=}~\mathcal{H}_y\,\mathfrak{D}_{c(M)}\vec_{\CR}(B\,\odot\,\Theta_{M}),$$
   where $\mathcal{H}_y$ is given by (\ref{eq210}).
$\blacksquare$
 For a better understanding of Lemma \ref{Lm31}, we consider the following example.
 \begin{exam}
     Consider 
     $A=\bmatrix{5 & 7 &3\\
     3 & 5 & 7\\
     7 & 3& 5}\in \CR_3$  {and} $M=\bmatrix{3 & 0& 9\\ 9 & 3 & 0\\ 0 &9 &3}\in \CR_3.$ 
     Then, $\Theta_M=\bmatrix{1 & 0 &1\\ 1 & 1& 0\\ 0&1&1},$ $\vec_{\CR}(\Theta_M)=\bmatrix{1\\1\\0},$ and 
     we get
     \begin{eqnarray}\label{exam:eq}
          {(A\,\odot\,\Theta_{M})x}=\left(\bmatrix{5 & 7 &3\\
     3 & 5 & 7\\
     7 & 3& 5}\odot \bmatrix{1 & 0 &1\\ 1 & 1& 0\\ 0&1&1} \right)\bmatrix{x_1\\x_2\\x_3}.
     \end{eqnarray}
     Then \eqref{exam:eq} can be rearranged in the following form:
     \begin{eqnarray}\label{exam:eq2}
       &&   {(A\,\odot\,\Theta_{M})x=} \bmatrix{x_1 & x_3 &x_2\\
     x_2 & x_1 & x_3\\
     x_3 & x_2& x_1} \bmatrix{1 & 0 &0\\ 0 & 1& 0\\ 0&0&0} \bmatrix{5\\3\\0}.
     \end{eqnarray}
    The above equation can be written as: $  {(A\,\odot\,\Theta_{M})x}= {\tt Cr}(x)\mathfrak{D}_{c(M)}\vec_{\CR}(A\odot \Theta_M).$
 \end{exam}
Next, we present the main result of this section concerning the structured \textit{BE} for  {\textit{circulant}} structured \textit{SPP} while preserving the sparsity pattern. Before that, we introduce the following notation:
 \begin{align}
&\mathfrak{D}_{c(A)}=\diag(\vec_{\CR}(\Theta_A)), \quad \mathfrak{D}_{c(B)}=\diag(\vec_{\CR}(\Theta_B)),\\
     &\mathfrak{D}_{c(C)}=\diag(\vec_{\CR}(\Theta_C)),\quad \mathfrak{D}_{\bm{a}}=\diag(\bm{a}),
 \end{align}
 where $\bm{a}=[\sqrt{n},\ldots,\sqrt{n}]^T\in \R^n.$
 
 \begin{theorem}\label{th31}
Let $\widetilde{u}=[\x^T,\, \y^T]^T$ be the approximate solution of the  {\textit{circulant}} structured \textit{SPP} \eqref{eq11}, i.e., $A,B,C\in \CR_n, $  and $w_4, w_5\neq 0.$ Then, we have 
\begin{align}\label{circulantBE}
    \eta^{\mathcal{S}_1}_{\tt sps}(\x,\y)=\left\|\mathcal{X}_{\tt Cr}^{H}(\mathcal{X}_{\tt Cr}\mathcal{X}_{\tt Cr}^H)^{-1} {r_{d}}\right\|_2,
\end{align}
where $\mathcal{X}_{\tt Cr}\in \C^{2n\times 5n}$ is given by
 \begin{align*}
    \mathcal{X}_{\tt Cr}=\bmatrix{ \frac{1}{w_1}{\tt Cr}\left(\x\right) \mathfrak{D}_{c(A)}\mathfrak{D}_{\bm{a}}^{-1} & \frac{1}{w_2}\mathcal{H}_\y\mathfrak{D}_{c(B)} \mathfrak{D}_{\bm{a}}^{-1} &  {\bf 0} &-\frac{1}{w_4}I_n &{\bf 0}\\{\bf 0}&\frac{1}{w_2} {\tt Cr}(\x) \mathfrak{D}_{c(B)}\mathfrak{D}_{\bm{a}}^{-1}& \frac{1}{w_3}{\tt Cr}\left(\y\right)\mathfrak{D}_{c(C)}\mathfrak{D}_{\bm{a}}^{-1} & {\bf 0} &-\frac{1}{w_5}I_n},
\end{align*}
$ {r_{d}}=\bmatrix{r_f^T, &r_g^T}^T,$ $r_f=f-A\x-B^T\y,$ and $r_g=g-B\x-C\y.$ 

Furthermore, the minimum norm perturbations to the  {Problem \ref{pr1}} are given by
\begin{align}\label{Eq32}
    &\widehat{\D A}_{\tt sps}={\tt Cr}\left(\frac{1}{w_1}\mathfrak{D}^{-1}_{\bm{a}}\bmatrix{I_n & {\bf 0}_{n\times 4n}}\mathcal{X}_{\tt Cr}^{H}(\mathcal{X}_{\tt Cr}\mathcal{X}_{\tt Cr}^H)^{-1} {r_{d}}\right),\\\label{Eq33}
    &\widehat{\D B}_{\tt sps}={\tt Cr}\left(\frac{1}{w_2}\mathfrak{D}^{-1}_{\bm{a}}\bmatrix{{\bf 0}_{n\times n}& I_n & {\bf 0}_{n\times 3n}}\mathcal{X}_{\tt Cr}^{H}(\mathcal{X}_{\tt Cr}\mathcal{X}_{\tt Cr}^H)^{-1} {r_{d}}\right),\\\label{Eq34}
    &\widehat{\D C}_{\tt sps}={\tt Cr}\left(\frac{1}{w_3}\mathfrak{D}^{-1}_{\bm{a}}\bmatrix{{\bf 0}_{n\times 2n}& I_n & {\bf 0}_{n\times 2n}}\mathcal{X}_{\tt Cr}^{H}(\mathcal{X}_{\tt Cr}\mathcal{X}_{\tt Cr}^H)^{-1} {r_{d}}\right),\\ \label{Eq35}
    &\widehat{\D f }_{\tt sps}=\frac{1}{w_4}\bmatrix{{\bf 0}_{n\times 3n}& I_n & {\bf 0}_{n\times n}}\mathcal{X}_{\tt Cr}^{H}(\mathcal{X}_{\tt Cr}\mathcal{X}_{\tt Cr}^H)^{-1} {r_{d}}, \quad \text{and}\\ \label{Eq36}
    &\widehat{\D g}_{\tt sps}=\frac{1}{w_5}\bmatrix{{\bf 0}_{n\times 4n}& I_n}\mathcal{X}_{\tt Cr}^{H}(\mathcal{X}_{\tt Cr}\mathcal{X}_{\tt Cr}^H)^{-1} {r_{d}}.
\end{align}
     \end{theorem}
     \proof Let $\widetilde{u}=[\widetilde{x}^T,\widetilde{y}^T]^T$ be an  approximate solution of the  {\textit{circulant}} structured \textit{SPP} of the form \eqref{eq11}. We need to construct sparsity preserving perturbations $
       \D A,\D B,
       \D C\in \CR_n,$ and perturbations $ \D f\in \C^n$ and $ \D g\in \C^n
    $. By Definition \ref{def24}, 
   $\D A, \D B,\D C, \D f,$ and $\D g$ satisfy 
     \begin{align}\label{eq31}
         &\D A\widetilde{x}+\D B^T \y-\D f=r_f,\\ \label{eq32}
         &\D B\widetilde{x}+\D C \y-\D g=r_g.
     \end{align}
     To maintain the sparsity pattern of $A,B$ and $C$ on the perturbation matrices, we replace  $\D A,\D B$ and
     $\D C $  by $\D A\bdot \Theta_A,\D B\bdot \Theta_B$ and $\D C\bdot \Theta_C,$ respectively.
     Consequently, from (\ref{eq31}) and (\ref{eq32}), we get 
     \begin{align}\label{eq33}
         &w_1^{-1}w_1(\D A\bdot \Theta_A)\widetilde{x}+w_2^{-1}w_2(\D B\bdot \Theta_B)^T \y-w_4^{-1}w_4\D f=r_f,\\ \label{eq34}
         &w_2^{-1}w_2(\D B\bdot \Theta_B)\widetilde{x}+w_3^{-1}w_3(\D C\bdot \Theta_C) \y-w_5^{-1}w_5\D g=r_g.
     \end{align}
      Applying Lemma \ref{Lm31} in  \eqref{eq33}, we obtain
      \begin{align}\label{eq35}
       \nonumber w_1^{-1} {\tt Cr}&\left(\x \right)\mathfrak{D}_{c(A)}w_1\vec_{\CR}(\D A\,\odot\,\Theta_{A})\\
       &+w_2^{-1}\mathcal{H}_\y\mathfrak{D}_{c(B)}w_2\vec_{\CR}(\D B\,\odot\,\Theta_{B})-w_4^{-1}w_4\D f=r_f.
      \end{align}
      Multiplying and dividing by $\mathfrak{D}_{\bm{a}}$ in \eqref{eq35}, we get 
      \begin{align}\label{eq36}
       \nonumber w_1^{-1} {\tt Cr}&\left(\x \right)\mathfrak{D}_{c(A)}\mathfrak{D}_{\bm{a}}^{-1}\mathfrak{D}_{\bm{a}}w_1\vec_{\CR}(\D A\,\odot\,\Theta_{A})\\
       &+w_2^{-1}\mathcal{H}_\y\mathfrak{D}_{c(B)}\mathfrak{D}_{\bm{a}}^{-1}\mathfrak{D}_{\bm{a}}w_2\vec_{\CR}(\D B\,\odot\,\Theta_{B})-w_4^{-1}w_4\D f=r_f.
      \end{align}
      We can reformulate (\ref{eq36}) as follows:
      \begin{equation}\label{eq38}
          \mathcal{X}_1\D \E=r_f,
      \end{equation}
      where \begin{equation*}
          \mathcal{X}_1=\bmatrix{ w_1^{-1}{\tt Cr}\left(\x \right)\mathfrak{D}_{c(A)}\mathfrak{D}_{\bm{a}}^{-1} &  w_2^{-1}\mathcal{H}_\y\mathfrak{D}_{c(B)}\mathfrak{D}_{\bm{a}}^{-1} &  {\bf 0} &-w_4^{-1}I_n &{\bf 0}}\in \C^{n\times 5n}
      \end{equation*}
      and \begin{equation}\label{EQ313}
          \D \E=\bmatrix{w_1\mathfrak{D}_{\bm{a}}\vec_{\CR}(\D A\,\odot\,\Theta_{A})\\w_2\mathfrak{D}_{\bm{a}}\vec_{\CR}(\D B\,\odot\,\Theta_{B})\\w_3\mathfrak{D}_{\bm{a}}\vec_{\CR}(\D C\,\odot\,\Theta_{C})\\ w_4\D f\\ w_5\D g}\in \C^{5n}.
      \end{equation}
      Here, the matrix $\mathfrak{D}_{\a}$ satisfy $\|\mathfrak{D}_{\a}\vec_{\CR}(A)\|_2=\|A\|_F,$ for any $A\in \CR_n.$
      
      Similarly, applying Lemma \ref{Lm31} to (\ref{eq34}), we obtain
       \begin{align}\label{eq37}
        \nonumber  w_2^{-1}{\tt Cr }&(\x ) \mathfrak{D}_{c(B)}w_2\vec_{\CR}(\D B\,\odot\,\Theta_{B})\\
          &+w_3^{-1}{\tt 
 Cr}\left(\y\right)\mathfrak{D}_{c(C)}w_3\vec_{\CR}(\D C\,\odot\,\Theta_{C})-w_5^{-1}w_5\D g=r_g.
      \end{align}
      Thus, we can reformulate (\ref{eq37}) as follows:
      \begin{equation}\label{eq310}
        \mathcal{X}_2\D\E=r_g,
      \end{equation}
      where
          $\mathcal{X}_2=\bmatrix{{\bf 0}& w_2^{-1}{\tt Cr}(\x) \mathfrak{D}_{c(B)}\mathfrak{D}_{\a}^{-1}&w_3^{-1}{\tt Cr}\left(\y \right)\mathfrak{D}_{c(C)}\mathfrak{D}_{\a}^{-1} & {\bf 0} &-w_5^{-1}I_n}\in \C^{n\times 5n}.$
      
\noindent      Now combining (\ref{eq38}) and (\ref{eq310}), we get the following equivalent linear system of (\ref{eq31})--(\ref{eq32}):
      \begin{equation}\label{eq311}
       \mathcal{X}_{\tt Cr}\D \E   \triangleq \bmatrix{\mathcal{X}_1\\ \mathcal{X}_2}\D \E= {r_{d}},
      \end{equation}
      where $\mathcal{X}_{\tt Cr}=\bmatrix{\mathcal{X}_1\\ \mathcal{X}_2}\in \C^{2n\times 5n}.$
Clearly, for $w_4,w_5\neq 0,$   the matrix $\mathcal{X}_{\tt Cr}$ has full row rank. Consequently, the linear system (\ref{eq311}) is consistent, and $\mathcal{X}^{\dagger}_{\tt Cr}=\mathcal{X}_{\tt Cr}^{H}(\mathcal{X}_{\tt Cr}\mathcal{X}_{\tt Cr}^H)^{-1},$ and by  {\cite[Theorem $1.1.6$]{GINVERSE}}, the minimum norm solution of the  system  is given by
      \begin{align}\label{eq312}
          \D \E_{\min}:=\mathcal{X}^{\dagger}_{\tt Cr} {r_{d}}=\mathcal{X}_{\tt Cr}^{H}(\mathcal{X}_{\tt Cr}\mathcal{X}_{\tt Cr}^H)^{-1} {r_{d}}.
      \end{align}
    Now, the minimization problem in Definition \ref{def24} is equivalently written as
    \begin{align}\label{eq313}
        \nonumber [\eta^{\mathcal{S}_1}_{\tt sps}(\widetilde{x},\widetilde{y})]^2&=\displaystyle{\min}\Big\{w_1^2\|\D A\bdot \Theta_A\|^2_F+w_2^2\|\D B\bdot \Theta_B\|_F^2 +w_3^3\|\D C\bdot \Theta_C\|_F^2  +w_4^2\|\D f\|_2^2+w_5^2\| \D g\|_2^2\, {\big |}\\
        &\hspace{3cm}{\bmatrix{
       \D A\bdot \Theta_A&(\D B\bdot \Theta_B)^T\\
    \D B\bdot \Theta_B &\D C\bdot \Theta_C}  
     \in \mathcal{S}_1}, \D f, \D g\in \C^n\Big\}\\
        &=\min \Bigg\{\|\D \E\|_2^2\,{\Big |}\mathcal{X}_{\tt Cr}\D\E= {r_{d}}\Bigg\}=\|\D \E_{\min}\|^2_2.
    \end{align}
     Hence, using (\ref{eq312}) and (\ref{eq313}), the structured \textit{BE} is given by 
     \begin{equation*}
         \eta^{\mathcal{S}_1}_{\tt sps}(\widetilde{x},\widetilde{y})=\left\|\mathcal{X}_{\tt Cr}^{H}(\mathcal{X}_{\tt Cr}\mathcal{X}_{\tt Cr}^H)^{-1} {r_{d}}\right\|_2.
     \end{equation*}
     From \eqref{EQ313}, we get $w_1\mathfrak{D}_{\bm{a}}\vec_{\CR}(\D A\bdot \Theta_{A})=\bmatrix{I_n & {\bf 0}_{n\times 4n}}\D \E.$ Thus, the minimal perturbation $\widehat{\D A}_{\tt sps}$ is given by
     $$\widehat{\D A}_{\tt sps}={\tt Cr}(\frac{1}{w_1}\mathfrak{D}^{-1}_{\bm{a}}\bmatrix{I_n & {\bf 0}_{n\times 4n}}\D \E_{\min}).$$
      Similarly, we can obtain other minimal perturbations.
     Hence, the proof is completed.
   $\blacksquare$

     In the following corollary, we present an explicit formula for the structured \textit{BE} $\eta^{\mathcal{S}_1}(\x,\y)$ for the  {\textit{circulant}} structured \textit{SPP} (\ref{eq11}) without maintaining the sparsity pattern in the perturbation matrices.
     \begin{corollary}
      Let $\widetilde{u}=[\x^T,\, \y^T]^T$ be the approximate solution of the  {\textit{circulant}} structured \textit{SPP}, i.e., $A,B,C\in \CR_n, $ and $w_4, w_5\neq 0.$ Then, we have
         \begin{equation}
             \eta^{\mathcal{S}_1}(\x,\y)=\left\|\widehat{\mathcal{X}}_{\tt Cr}^H(\widehat{\mathcal{X}}_{\tt Cr}\widehat{\mathcal{X}}_{\tt Cr}^H)^{-1} {r_{d}}\right\|_2,
         \end{equation}
         where \begin{equation}\label{eq320}
             \widehat{\mathcal{X}}_{\tt Cr}=\bmatrix{\frac{1}{w_1}{\tt Cr}\left(\x\right)\mathfrak{D}_{\bm{a}}^{-1} & \frac{1}{w_2}\mathcal{H}_\y\mathfrak{D}_{\bm{a}}^{-1} &  {\bf 0} &-\frac{1}{w_4}I_n &{\bf 0}\\{\bf 0}&\frac{1}{w_2} {\tt Cr}\left(\x\right)\mathfrak{D}_{\bm{a}}^{-1}& \frac{1}{w_3}{\tt Cr}\left(\y\right)\mathfrak{D}_{\bm{a}}^{-1} & {\bf 0} &-\frac{1}{w_5}I_n}.
         \end{equation}
     \end{corollary}
     \proof
         Since we are not maintaining the sparsity pattern to the perturbation matrices, we consider $\Theta_A=\Theta_B=\Theta_C={\bf 1}_{n\times n}.$ Then $\D A \bdot \Theta_A=\D A,$ $\D B \bdot \Theta_B=\D B$ and $\D C \bdot \Theta_C=\D C.$ Also, $\vec_{\CR}(\Theta_A)=\vec_{\CR}(\Theta_B)=\vec_{\CR}(\Theta_C)={\bf 1}_n.$ Consequently, the proof is completed using the formula stated in Theorem \ref{th31}.
       $\blacksquare$

         The minimal perturbations $\widehat{\D A}, \widehat{\D B}, \widehat{\D C}, \widehat{\D f}, $ and $\widehat{\D g}$ to the Problem \ref{pr1} are given by formulae (\ref{Eq32})-(\ref{Eq36}), where $\mathcal{X}_{\tt Cr}=\widehat{\mathcal{X}}_{\tt Cr}$.

     \section{Structured \textit{BEs} 
     for \textit{Toeplitz} structured \textit{SPPs}}\label{SEC4}
    This section focuses on the derivation of compact formulae for the structured \textit{BEs} $\eta^{\mathcal{S}_2}_{\tt sps}(\x,\y)$ and $\eta^{\mathcal{S}_2}(\x,\y)$ for \textit{Toeplitz} structured \textit{SPPs} with and without preserving sparsity pattern, respectively. 
    In addition, the minimal perturbations are provided for the Problem \ref{pr1} for which the structured \textit{BEs} are obtained. To accomplish this, we first derive the following lemma.
      \begin{lemma}\label{lm41}
     Let $A,B,M\in\T_{ m\times n}$ with  generator vectors $\vec_{\T}(A)=[a_{- {m}+1},\ldots,a_{-1},a_0,a_1\ldots,a_{n-1}]^T\\ {\in \C^{n+m-1}},$ $\vec_{\T}(B)=[b_{- {m}+1},\ldots,b_{-1},b_0,b_1\ldots,b_{n-1}]^T$$ {\in \C^{n+m-1}},$ and  $\vec_{\T}(M)=[m_{- {m}+1},\ldots,m_{-1},\\m_0,m_1\ldots,m_{n-1}]^T {\in \C^{n+m-1}},$ respectively. Suppose $x=[x_1,\ldots,x_n]^T\in \C^n,$  ${y=[y_1,\ldots,y_{ {m}}]^T\in\C^{ {m}}}.$ 
     Then 
     { $$(A\,\odot\,\Theta_{M})x=\mathcal{K}_x\mathfrak{D}_{t(M)}\vec_{\T}(A\,\odot\,\Theta_{M})\quad \text{and}$$ 
     $$(B\,\odot\,\Theta_{M})^Ty= \mathcal{G}_y\mathfrak{D}_{t(M)}\vec_{\T}(B\,\odot\,\Theta_{M}),$$}
     where $t(M)=\vec_{\T}(\Theta_M),$

\[
 \mathcal{K}_x = 
  \mathop{\left[
  \begin{array}{ *{10}{c} }
    0 & \cdots& \cdots& 0 & \colind{x_1}{\tiny $m^{th}$ term }& \cdots&\cdots&\cdots &x_{n-1}&x_n\\
                                \vdots&  \cdots &0& x_1& x_2&\cdots& \cdots&x_{n-1}&x_n&0\\
                                \vdots&\iddots &\iddots&\iddots & & &\iddots &\iddots&\iddots&\vdots\\
                               0 & \iddots&\iddots & & &\iddots &\iddots&\iddots &&\vdots\\
                                x_1&x_2&\cdots& &x_{n-1}&x_n&0&\cdots&\cdots&0
  \end{array}
  \right]}
  \in \C^{ m\times (m+n-1)},
\]
\[
 \mathcal{G}_y =\color{black}
  \mathop{\left[
  \begin{array}{ *{9}{c} }
   y_m & y_{m-1}&\cdots&y_1&0&\cdots&\cdots&\cdots&0\\
                                0&y_m&y_{m-1}&\cdots &y_1&0&\cdots&\cdots&0\\
                                0& 0&\ddots & \ddots& \ddots&\ddots & & & \vdots\\
                                \vdots& & \ddots & \ddots&\ddots & \ddots&\ddots & &  \vdots\\
                               0 &\cdots &  & 0&y_m &y_{m-1} & \cdots& y_1&0  \\
                                0&\cdots&\cdots&\cdots&0&\rowind{y_m }{\tiny $n^{th}$ term}&y_{m-1}&\cdots&y_1
  \end{array}
  \right]}
  \in \C^{n\times (m+n-1)}.
\]
    
 \end{lemma}
  \Proof
    The proof proceeds in a similar manner to the proof of Lemma \ref{Lm31}.
$\blacksquare$

In the following theorem, we derive the explicit formula for the structured \textit{BE} $\eta^{\mathcal{S}_2}_{\tt sps}(\x,\y)$. 
Prior to that, we introduce the following notation:
 \begin{align}
&\mathfrak{D}_{t(A)}=\diag(\vec_{\T}(\Theta_A)), \quad \mathfrak{D}_{t(B)}=\diag(\vec_{\T}(\Theta_B)),\\
     &\mathfrak{D}_{t(C)}=\diag(\vec_{\T}(\Theta_C))\quad \text{and}\quad  {\mathfrak{D}_{{\bf t}_{mn}}=\diag({\bf t}_{mn})},
 \end{align}
 where $ {{\bf t}_{mn}}=[1,\sqrt{2},\ldots, {\sqrt{m-1}}, {\sqrt{\min\{m,n\}}},\sqrt{n-1},\ldots,\sqrt{2},1]^T\in \R^{ m+n-1}.$  {When $n=m,$ we write ${\bf t}_{mn}={\bf t}_{n}$ (or ${\bf t}_{m}$)}. 
 \begin{theorem}\label{th41}
     Let $\widetilde{u}=[\x^T,\, \y^T]^T$ be the approximate solution of the \textit{Toeplitz} structured \textit{SPP} \eqref{eq11}, i.e.,  {$A\in \T_{n\times n}, B\in \T_{m\times n}, C\in \T_{m\times m} $} and $w_4, w_5\neq 0.$ Then, we have 
\begin{align}
    \eta^{\mathcal{S}_2}_{\tt sps}(\x,\y)=\left\|\mathcal{X}_{\T}^{H}(\mathcal{X}_{\T}\mathcal{X}_{\T}^H)^{-1} {r_{d}}\right\|_2,
\end{align}
where $ \mathcal{X}_{\T}\in \C^{ (n+m)\times (4n+4m-3)}$ is given by
 \begin{align*}
    \mathcal{X}_{\T}=\bmatrix{ \frac{1}{w_1}\mathcal{K}_\x\mathfrak{D}_{t(A)} {\mathfrak{D}_{\t_n}^{-1}} & \frac{1}{w_2}\mathcal{G}_\y\mathfrak{D}_{t(B)} {\mathfrak{D}_{\t_{mn}}^{-1}} &  {\bf 0} &-\frac{1}{w_4}I_n &{\bf 0}\\{\bf 0}&\frac{1}{w_2} \mathcal{K}_\x\mathfrak{D}_{t(B)} {\mathfrak{D}_{\t_{mn}}^{-1}}& \frac{1}{w_3}\mathcal{K}_\y\mathfrak{D}_{t(C)} {\mathfrak{D}_{\t_m}^{-1}} & {\bf 0} &-\frac{1}{w_5}I_m},
\end{align*}
$ {r_{d}}=\bmatrix{r_f^T, &r_g^T}^T,$ $r_f=f-A\x-B^T\y,$ and $r_g=g-B\x-C\y.$

Furthermore, the minimal  perturbations to the Problem \ref{pr1} are given by
{\footnotesize\begin{align}\label{Eq42}
    &\widehat{\D A}_{\tt sps}={\T}\left(\frac{1}{w_1} {\mathfrak{D}^{-1}_{\t_n}}\bmatrix{I_{2n-1} & {\bf 0}_{(2n-1)\times  (2n+4m-2)}}\mathcal{X}_{\T}^{H}(\mathcal{X}_{\T}\mathcal{X}_{\T}^H)^{-1} {r_{d}}\right),\\\label{Eq43}
    & \widehat{\D B}_{\tt sps}={\T}\left(\frac{1}{w_2} {\mathfrak{D}^{-1}_{\t_{mn}}}\bmatrix{{\bf 0}_{ (m+n-1)\times (2n-1)}& I_{ m+n-1} & {\bf 0}_{ (m+n-1)\times (3m+n-1)}}\mathcal{X}_{\T}^{H}(\mathcal{X}_{\T}\mathcal{X}_{\T}^H)^{-1} {r_{d}}\right),\\\label{Eq44}
    &\widehat{\D C}_{\tt sps}={\T}\left(\frac{1}{w_3} {\mathfrak{D}^{-1}_{\t_m}}\bmatrix{{\bf 0}_{ (2m-1)\times (3n+3m-2)}& I_{ 2m-1} & {\bf 0}_{ (2m-1)\times (n+m)}}\mathcal{X}_{\T}^{H}(\mathcal{X}_{\T}\mathcal{X}_{\T}^H)^{-1} {r_{d}}\right),\\ \label{Eq45}
    &\widehat{\D f }_{\tt sps}=\frac{1}{w_4}\bmatrix{{\bf 0}_{  n\times (3n+3m-3)}& I_n & {\bf 0}_{ n\times m}}\mathcal{X}_{\T}^{H}(\mathcal{X}_{\T}\mathcal{X}_{\T}^H)^{-1} {r_{d}},\, \text{and}\\ \label{Eq46}
   & \widehat{\D g}_{\tt sps}=\frac{1}{w_5}\bmatrix{{\bf 0}_{ m\times (4n+3m-3)}& I_{ m}}\mathcal{X}_{\T}^{H}(\mathcal{X}_{\T}\mathcal{X}_{\T}^H)^{-1} {r_{d}}.
\end{align}}
 \end{theorem}
 \proof
 We need to construct perturbations  {$\D A\in \T_{n\times n},\D B\in \T_{m\times n},
     \D C\in \T_{m\times m}$} (which preserves the sparsity pattern of the original matrices), $\D f\in \C^n ~\text{and}~ \D g\in \C^{ {m}}
     $ for the approximate solution $\widetilde{u}=[\widetilde{x}^T,\widetilde{y}^T]^T.$ By   {Definition \ref{def24},} 
     $\D A, \D B,\D C, \D f,$ and $\D g$ satisfy 
     \begin{align*}
         &\D A\widetilde{x}+\D B^T \y-\D f=r_f,\\ \label{eq42}
         &\D B\widetilde{x}+\D C \y-\D g=r_g.
     \end{align*}
     Following the proof method of Theorem \ref{th31} and using Lemma \ref{lm41}, we obtain
     { $ {r_{d}}=\mathcal{X}_{\T}\D \E,$
     where }  \begin{equation}\label{Tp_eq413}
          \D \E=\bmatrix{w_1 {\mathfrak{D}_{\t_n}}\vec_{\T}(\D A\,\odot\,\Theta_{A})\\w_2 {\mathfrak{D}_{\t_{mn}}}\vec_{\T}(\D B\,\odot\,\Theta_{B})\\w_3 {\mathfrak{D}_{\t_m}}\vec_{\T}(\D C\,\odot\,\Theta_{C})\\ w_4\D f\\ w_5\D g}\in \C^{ 4n+4m-3}.
      \end{equation} 
     Hence, an analogous way to proof of Theorem \ref{th31}, we get $$\eta_{\tt sps}^{\mathcal{S}_2}(\x,\y)=\left\|\mathcal{X}_{\T}^{H}(\mathcal{X}_{\T}\mathcal{X}_{\T}^H)^{-1} {r_{d}}\right\|_2.$$
      From \eqref{Tp_eq413}, we get $w_1 {\mathfrak{D}_{\t_n}}\vec_{\T}(\D A\,\odot\,\Theta_{A})=\bmatrix{I_{2n-1}& {\bf 0}_{{(2n-1)}\times { (2n+4m-2)}}}\D \E.$ Therefore,  the minimal perturbation matrix $\widehat{\D A}_{\tt sps}$ which preserve the sparsity pattern of $A$ is given by
      $$\widehat{\D A}_{\tt sps}=\frac{1}{w_1} {\mathfrak{D}_{\t_n}^{-1}}\bmatrix{I_{2n-1}& {\bf 0}_{{(2n-1)}\times { (2n+4m-2)}}}\D \E_{\min}.$$
      Similarly, we can obtain the minimal perturbations $\widehat{\D B}_{\tt sps}, \widehat{\D C}_{\tt sps}, \widehat{\D f }_{\tt sps}$ and $\widehat{\D g}_{\tt sps}.$ Hence, the proof is completed. $\blacksquare$

  The next corollary presents a formula for $\eta^{\mathcal{S}_2}(\x,\y)$ without considering the sparsity pattern in the input matrices.
 \begin{corollary}\label{coro41}
 Let $\widetilde{u}=[\x^T,\, \y^T]^T$ be the approximate solution of the \textit{Toeplitz}  structured \textit{SPP} \eqref{eq11}, i.e., $A,B,C\in \T_n,$ and $w_4, w_5\neq 0.$  Then, we have 
     \begin{equation}
             \eta^{\mathcal{S}_2}(\x,\y)=\left\|\widehat{\mathcal{X}}_{\T}^{H}(\widehat{\mathcal{X}}_{\T}\widehat{\mathcal{X}}_{\T}^H)^{-1} {r_{d}}\right\|_2,
         \end{equation}
         where $ \widehat{\mathcal{X}}_{\T}\in \C^{ (n+m)\times (4n+4m-3)}$ is given by\begin{align}\label{eq420}
            \widehat{\mathcal{X}}_{\T}=\bmatrix{\frac{1}{w_1}\mathcal{K}_{\x} {\mathfrak{D}_{\t_n
             }^{-1}} & \frac{1}{w_2}\mathcal{G}_\y {\mathfrak{D}_{\t_{mn}}^{-1} }& {\bf 0}&-\frac{1}{w_4}I_n &{\bf 0}\\{\bf 0}&\frac{1}{w_2} \mathcal{K}_{\x} {\mathfrak{D}_{\t_{mn}}^{-1}}& \frac{1}{w_3}\mathcal{K}_{\y} {\mathfrak{D}_{\t_m}^{-1}} & {\bf 0} &-\frac{1}{w_5}I_m}.
         \end{align}
 \end{corollary}
 \proof
Since the sparsity pattern of the perturbation matrices is not taken care of, we consider  {$\Theta_A={\bf 1}_{n\times n},$ $ \Theta_B={\bf 1}_{m\times n} $ and  $\Theta_C={\bf 1}_{m\times m}.$} Then $\D A \bdot \Theta_A=\D A,$ $\D B \bdot \Theta_B=\D B,$ and $\D C \bdot \Theta_C=\D C.$ Also,  {$\vec_{\T}(\Theta_A)={\bf 1}_{2n-1},$ $ \vec_{\T}(\Theta_B)={\bf 1}_{m+n-1}$ and  $ \vec_{\T}(\Theta_C)={\bf 1}_{2m-1}.$} As a result, the proof follows using the formula stated in Theorem \ref{th41}.     
$\blacksquare$

  { Note that the} minimal perturbations $\widehat{\D A}, \widehat{\D B}, \widehat{\D C}, \widehat{\D f}, $ and $\widehat{\D g}$ to the Problem \ref{pr1} are given by formulae (\ref{Eq42})-(\ref{Eq46}) with $\mathcal{X}_{\T}=\widehat{\mathcal{X}}_{\T}.$     
 \section{Structured \textit{BEs} 
 for \textit{symmetric}-\textit{Toeplitz} structured \textit{SPPs}}\label{SEC5}
In this section, we derive concise formulae for the structured  \textit{BE} $\eta^{\mathcal{S}_3}_{\tt sps}(\x,\y)$ and $\eta^{\mathcal{S}_3}(\x,\y)$  for \textit{symmetric}-\textit{Toeplitz} structured \textit{SPPs} with and without preserving the sparsity pattern, respectively.  {Since the \( (1,2) \) block matrix \( B \) is symmetric, thus the case \( n = m \) follows.} 
Additionally, we derive the structured \textit{BE} when block matrices $A, B$, and $C$ only preserve the sparsity pattern.
 \subsection{\textit{Symmetric}-\textit{Toeplitz} structured \textit{SPPs}}
 In many applications, such as the \textit{WRLS} problem, the block matrices $A$ and $C$ do not follow any particular structure, in this subsection, we focus on the structured \textit{BE} when the perturbation matrix $\D B$ follow \textit{symmetric}-\textit{Toeplitz} structure of $B$.  To find the structured $BE,$ 
 we present the following lemmas that are crucial in establishing our main results. 
 \begin{lemma}\label{lm51}
     Let $A,M\in {\mathcal{S}\T}_n$ with  generator vectors $\vec_{\mathcal{S}\T}(A)=[a_0,a_1\ldots,a_{n-1}]^T\in \C^n$  and $\vec_{\mathcal{S}\T}(M)=[m_0,m_1\ldots,m_{n-1}]^T\in \C^n,$ respectively. Suppose $x=[x_1,\ldots,x_n]^T\in \C^n,$  {then 
     $$(A\,\odot\,\Theta_{M})x=\mathcal{I}_x\mathfrak{D}_{s(M)}\vec_{\mathcal{S}\T}(A\,\odot\,\Theta_{M}),$$} 
     where $s(M)=\vec_{\mathcal{S}\T}(\Theta_M)$  and $\mathcal{I}_{x}\in \C^{n\times n}$ is given by \begin{align}
         \mathcal{I}_{x}=\bmatrix{x_1& \cdots&\cdots&x_{n-1}&x_n\\
                               x_2&\cdots&\cdots &x_n&0\\
                               \vdots&&\iddots&\iddots&\vdots\\
                               \vdots&x_n&0&\cdots&\vdots\\
                               x_{n} &0&\cdots&\cdots&0
         }+\bmatrix{0 &\cdots&\cdots&&0\\
                    0&x_1&0&\cdots&0\\
                    \vdots&\vdots&\ddots&\ddots&\vdots\\
                    \vdots&x_{n-2}&\cdots&x_1&0\\
                    0&x_{n-1}&\cdots&\cdots&x_1
                            }.
     \end{align}
 \end{lemma}
 \Proof
    The proof proceeds in a similar manner to the proof of Lemma \ref{Lm31}.
$\blacksquare$

 \begin{lemma}\label{lm52}
       Let $A,B,M\in \C^{m\times n}$ be three matrices. Suppose that $x=[x_1,\ldots,x_n]^T\in \C^n$ and $y=[y_1,\ldots,y_m]^T\in \C^m.$  
       Then 
       \begin{eqnarray*}
           { (A\,\odot\,\Theta_{M})x=\mathcal{J}^{1,m}_x\mathfrak{D}_{\vec(\Theta_M)}\vec(A\,\odot\,\Theta_{M})\, \mbox{and }\, (B\bdot \Theta_M)^Ty=\mathcal{J}^{2,n}_{y}\mathfrak{D}_{\vec(\Theta_M)}\vec(B\,\odot\,\Theta_{M}),}
       \end{eqnarray*}
  \noindent where 
    {$\mathcal{J}^{1,m}_x=x^T\otimes I_m\in \C^{m\times {mn}},$  $\mathcal{J}^{2,n}_{y}=I_n\otimes y^T\in \C^{n\times {mn}},$
   and $\otimes$ denotes the Kronecker product.}
 \end{lemma}
 \Proof
    The proof proceeds in a similar method to the proof of Lemma \ref{Lm31}.
$\blacksquare$

\vspace{1mm}
 Next, we derive concrete formulae for $\eta_{\tt sps}^{\mathcal{S}_3}(\x,\y)$ and $\eta^{\mathcal{S}_3}(\x,\y),$ which are the main result of this section.  Before proceeding, we introduce the following notations:
 \begin{align}
&\mathfrak{D}_{s(A)}=\diag(\vec_{\ST}(\Theta_A)), \quad \mathfrak{D}_{s(B)}=\diag(\vec_{\ST}(\Theta_b)),\\
     &\mathfrak{D}_{s(C)}=\diag(\vec_{\ST}(\Theta_C))\quad \text{and}\quad \mathfrak{D}_{{\bf s}}=\diag({\bf s}),
 \end{align}
 where ${\bf s}=[\sqrt{n},\sqrt{2(n-1)},\sqrt{2(n-2)},\ldots,\sqrt{2}]^T\in \R^n.$
 \begin{theorem}\label{th51}
     Let $\widetilde{u}=[\x^T,\, \y^T]^T$ be  an approximate solution of the \textit{symmetric}-\textit{Toeplitz} structured \textit{SPP} \eqref{eq11}, i.e., $B\in \ST_n,$  $A,C \in\C^{n\times n},$  and $w_4, w_5\neq 0.$  Then, we have 
\begin{align}
    \eta^{\mathcal{S}_3}_{\tt sps}(\x,\y)=\left\|\mathcal{Z}_{\ST}^{H}(\mathcal{Z}_{\ST}\mathcal{Z}_{\ST}^H)^{-1} {r_{d}}\right\|_2,
\end{align}
where $ \mathcal{Z}_{\ST}\in \C^{2n\times l}$ is given by
\begin{align}
 \mathcal{Z}_{\ST}=\bmatrix{\frac{1}{w_1} \mathcal{J}^{1,n}_{\x}\mathfrak{D}_{\vec(\Theta_A)}&\frac{1}{w_2} \mathcal{I}_{\y}\mathfrak{D}_{s(B)}\mathfrak{D}^{-1}_{\bf s}&{\bf 0} &-\frac{1}{w_4}I_n & {\bf 0}\\
 {\bf 0} &\frac{1}{w_2} \mathcal{I}_{\x}\mathfrak{D}_{s(B)}\mathfrak{D}^{-1}_{\bf s}& \frac{1}{w_3}\mathcal{J}^{1,n}_{\y}\mathfrak{D}_{\vec(\Theta_C)}&{\bf 0}& -\frac{1}{w_5}I_n},
\end{align}
 $ {r_{d}}=\bmatrix{r_f^T, &r_g^T}^T,$ $r_f=f-A\x-B\y,$ $r_g=g-B\x-C\y,$ and $l=2n^2+3n.$

Furthermore, the minimal perturbations to the Problem \ref{pr1} are given by
\begin{align}\label{Eq55}
    &\vec(\widehat{\D A}_{\tt sps})=\frac{1}{w_1}\bmatrix{I_{n^2} & {\bf 0}_{n^2\times (n^2+3n)}}\mathcal{Z}_{\ST}^{H}(\mathcal{Z}_{\ST}\mathcal{Z}_{\ST}^H)^{-1} {r_{d}},\\ \label{Eq56}
    &\widehat{\D B}_{\tt sps}={\ST}\left(\frac{1}{w_2}\mathfrak{D}^{-1}_{\bf s}\bmatrix{{\bf 0}_{n\times n^2}& I_{n} & {\bf 0}_{n\times (2n^2+2n)}}\mathcal{Z}_{\ST}^{H}(\mathcal{Z}_{\ST}\mathcal{Z}_{\ST}^H)^{-1} {r_{d}}\right),\\\label{Eq57}
    &\vec(\widehat{\D C}_{\tt sps})=\frac{1}{w_3}\bmatrix{{\bf 0}_{n^2\times (n^2+n)}& I_{n^2} & {\bf 0}_{n^2\times (n^2+2n)}}\mathcal{Z}_{\ST}^{H}(\mathcal{Z}_{\ST}\mathcal{Z}_{\ST}^H)^{-1} {r_{d}},\\
    &\widehat{\D f }_{\tt sps}=\frac{1}{w_4}\bmatrix{{\bf 0}_{n\times (2n^2+n)}& I_n & {\bf 0}_{n\times (2n^2+n)}}\mathcal{Z}_{\ST}^{H}(\mathcal{Z}_{\ST}\mathcal{Z}_{\ST}^H)^{-1} {r_{d}},\quad \text{and}\\\label{Eq58}
    &\widehat{\D g}_{\tt sps}=\frac{1}{w_5}\bmatrix{{\bf 0}_{n\times (2n^2+2n)}& I_n}\mathcal{Z}_{\ST}^{H}(\mathcal{Z}_{\ST}\mathcal{Z}_{\ST}^H)^{-1} {r_{d}}.
 \end{align}
 \end{theorem}
 \proof
For an approximate solution $\widetilde{u}=[\widetilde{x}^T,\widetilde{y}^T]^T,$ we require to construct sparsity preserving perturbations $
       \D B\in \ST_n, \D A, \D C\in \C^{n\times n},$ and perturbations $ \D f\in \C^n$ and $ \D g\in \C^n
    $. By  {Definition \ref{def24}, we have}
     \begin{align}\label{eq59}
         &\D A\widetilde{x}+\D B^T \y-\D f=r_f,\\ \label{eq510}
         &\D B\widetilde{x}+\D C \y-\D g=r_g,
     \end{align}
     where $\D B\in \ST_n, \D A,\D C\in \C^{n\times n}.$ \\    
     Following the proof method of Theorem \ref{th31} and using Lemmas \ref{lm51} and \ref{lm52}, we obtain  
      { $ {r_{d}}=\mathcal{Z}_{\ST}\D \E,$ where}
      \begin{equation}
          \D \E=\bmatrix{w_1\mathfrak{D}_{\bf s}\vec(\D A\,\odot\,\Theta_{A})\\w_2\mathfrak{D}_{\bf s}\vec_{\ST}(\D B\,\odot\,\Theta_{B})\\w_3\mathfrak{D}_{\bf s}\vec(\D C\,\odot\,\Theta_{C})\\ w_4\D f\\ w_5\D g}\in \C^{l}.
      \end{equation}

       {Hence, following an analogous way to proof of Theorem \ref{th31}, we get 
      the desired structured \textit{BE} and perturbation matrices.}
$\blacksquare$

 Next, we present the formula for  $\eta^{\mathcal{S}_3}(\x,\y)$ without preserving the sparsity pattern.
 \begin{corollary}\label{coro51}
     Let $\widetilde{u}=[\x^T,\, \y^T]^T$ be an approximate solution of the \textit{symmetric}-\textit{Toeplitz} structured \textit{SPP} \eqref{eq11}, i.e., $B\in \ST_n$  $A,C \in\C^{n\times n},$ and $w_4, w_5\neq 0.$  Then, we have 
     \begin{align}
             \eta^{\mathcal{S}_3}(\x,\y)=\left\|\widehat{\mathcal{Z}}_{\ST}^{H}(\widehat{\mathcal{Z}}_{\ST}\widehat{\mathcal{Z}}_{\ST}^H)^{-1} {r_{d}}\right\|_2,
             \end{align}
             where \begin{equation}
                   {\widehat{\mathcal{Z}}_{\ST}}=\bmatrix{\frac{1}{w_1}\mathcal{J}^{1,n}_{\x}&\frac{1}{w_2}\mathcal{I}_{\y}\mathfrak{D}^{-1}_{\bf s}&{\bf 0} &-\frac{1}{w_4}I_n & {\bf 0}\\
 {\bf 0} &\frac{1}{w_2}\mathcal{I}_{\x}\mathfrak{D}^{-1}_{\bf s}& \frac{1}{w_3}\mathcal{J}^{1,n}_{\y}&{\bf 0}& -\frac{1}{w_5}I_n}.
             \end{equation}
 \end{corollary}
 \proof
   Because we are not considering the sparsity pattern in the perturbation matrices, taking $\Theta_A=\Theta_B=\Theta_C={\bf 1}_{n\times n}$ in the Theorem \ref{th51} yields the desired result.
 $\blacksquare$
 
      The minimal perturbation matrices $\widehat{\D A}, \widehat{\D B}, \widehat{\D C}, \widehat{\D f}, $ and $\widehat{\D g}$ to the Problem \ref{pr1} are given by formulae (\ref{Eq55})-(\ref{Eq58}) with $\mathcal{Z}_{\ST}=\widehat{\mathcal{Z}}_{\ST}.$  
 \begin{remark}
     Applying our framework developed in this subsection and Sections \ref{SEC3} and \ref{SEC4}, we can obtain the structured \textit{BEs} for the SPP \eqref{eq11} when the block matrices possess only \textit{symmetric} structure or   {Hankel structure} (which is \textit{symmetric} as well).  
 \end{remark}
 \subsection{Unstructured \textit{BEs} with preserving sparsity pattern}
 In this subsection, we address the scenario where  $A$, $B$, and $C$ in \eqref{eq11} are unstructured. Although the previous studies such as \cite{be2007wei, be2012LAA,  BE2020BING} have explored the  \textit{BEs} for the \textit{SPP} \eqref{eq11}, their investigations do not take into account the sparsity pattern of the block matrices. Consequently, in this scenario, first we define the unstructured \textit{BE} as follows:
 \begin{equation}\label{eq516}
  \eta(\widetilde{x},\widetilde{y}):=\displaystyle{\min_{\substack{\D \mathcal{A}\in\, \mathcal{S}^0, \, \D d\in \C^{n+m},\\ (\A+\D \A)\widetilde{u}=d+\D d}}} \vertiii{\bmatrix{\D \A& \D{d}}}_{{ w},F}
 \end{equation}
 where  {$\mathcal{S}^0=\left\{ \D \A=\bmatrix{\D A& \D B^T\\ \D B & \D C} ~{\Big|}~ \D A\in \C^{n\times n}, \D B\in \C^{m\times n}, \D C\in \C^{m\times m}\right\}.$}

      {The following result gives} the formula for the unstructured \textit{BE} of the \textit{SPPs} $(\ref{eq11}) $ when the sparsity pattern in the perturbation matrices is preserved, and in this case, we denote it by $\eta_{\tt sps}(\x,\y)$.
     
     \begin{theorem}\label{th52}
     Let $\widetilde{u}=[\x^T,\, \y^T]^T$ be  an approximate solution of the \textit{SPP} (\ref{eq11}) with $A\in \C^{n\times n},B\in \C^{m\times n},$ $C\in \C^{m\times m},$ and $w_4, w_5\neq 0.$ Then, we have
     \begin{align*}
    \eta_{\tt sps}(\x,\y)=\left\|\mathcal{N}^{H}(\mathcal{N}\mathcal{N}^H)^{-1} {r_{d}}\right\|_2,
\end{align*}
where $ \mathcal{N}\in \C^{(m+n)\times k}$ is given by
\begin{align}
 \mathcal{N}=\bmatrix{\frac{1}{w_1}\mathcal{J}^{1,n}_{\x}\mathfrak{D}_{\vec(\Theta_A)}&\frac{1}{w_2}\mathcal{J}^{2,n}_{\y}\mathfrak{D}_{\vec(\Theta_B)}&{\bf 0} &-\frac{1}{w_4}I_n & {\bf 0}\\
 {\bf 0} &\frac{1}{w_2}\mathcal{J}^{1,m}_{\x}\mathfrak{D}_{\vec(\Theta_B)}&\frac{1}{w_3}\mathcal{J}^{1,m}_{\y}\mathfrak{D}_{\vec(\Theta_C)}&{\bf 0}& -\frac{1}{w_5}I_m}
\end{align}
 {and} $k=n^2+m^2+mn+n+m.$
\end{theorem}
\proof The proof follows similarly to the proof method of  Theorem \ref{th51}. $\blacksquare$


The next corollary presents the \textit{BE} formula when the sparsity is not considered.
\begin{corollary}\label{coro52}
Let $\widetilde{u}=[\x^T,\, \y^T]^T$ be  an approximate solution of the \textit{SPP} (\ref{eq11}) with $A\in \C^{n\times n},B\in \C^{m\times n}$ and $C\in \C^{m\times m},$ and $w_4, w_5\neq 0.$ Then, we have
\begin{align*}
    \eta(\x,\y)=\left\|\widehat{\mathcal{N}}^{H}(\widehat{\mathcal{N}}\widehat{\mathcal{N}}^H)^{-1} {r_{d}}\right\|_2,
\end{align*}
where $ \widehat{\mathcal{N}}\in \C^{(m+n)\times l}$ is given by
    \begin{align}
        \widehat{\mathcal{N}}=\bmatrix{\frac{1}{w_1}\mathcal{J}^{1,n}_{\x}&\frac{1}{w_2}\mathcal{J}^{2,n}_{\y}&{\bf 0} &-\frac{1}{w_4}I_n & {\bf 0}\\
 {\bf 0} &\frac{1}{w_2}\mathcal{J}^{1,m}_{\x}&\frac{1}{w_3}\mathcal{J}^{1,m}_{\y}&{\bf 0}& -\frac{1}{w_5}I_m}.
    \end{align}
\end{corollary}
\proof The proof is followed by taking $\Theta_A={\bf 1}_{m\times m},\, \Theta_B={\bf 1}_{m\times n}$ and $\Theta_C={\bf 1}_{n\times n}$ in the expression of $\eta_{\tt sps}(\x,\y),$  presented in Theorem \ref{th52}. $\blacksquare$

Note that when \( w_4 \) or \( w_5 \) are zero, the desired BE is achieved if \( \mathcal{N} \) and \( \widehat{\mathcal{N}} \) have full row rank. Nevertheless, in \cite{BEKKT2004, be2007wei}, formulas for \textit{BEs} with no special structure on block matrices are discussed,  the following example illustrates that our \textit{BE} can be smaller than theirs.
\begin{exam}
    Consider the \textit{SPP} \eqref{eq11}, where $A=I_4,$ $B=\bmatrix{2&1&3&1\\	-1 &	2&	1&	1}\in \R^{2\times 4},$ $C={\bf 0},$ $f=[-1, 0,2,3]^T,$ and $g={\bf 0}.$ We take the approximate solution $[\x^T,\, \y^T]^T=[-1.495, 1.505,1.505,1.505,1.005,-0.495]^T.$ Then, employing the formula provided in  \cite{BEKKT2004}  {with $\theta=1$}, the computed \textit{BE} is $0.0410$.  {Since, the (1,1) block in \cite{BEKKT2004} has no perturbation, by considering $w_1=0$  in Corollary \ref{coro52} (with $w_2=w_{ {4}}=1,$ $w_3=w_5=0$), the \textit{BE} is $0.0288.$ This comparison highlights that our computed \textit{BE} using Corollary \ref{coro52} and  \cite{BEKKT2004} are of the same order, illustrating the reliability of our obtained \textit{BE}.} 
    \end{exam}
\section{Application to derive the structured \textit{BEs} for the \textit{WRLS} problems}\label{SEC6}
In this section, we present  {an application} of our developed theory in deriving the structured \textit{BE} for the \textit{WRLS} problem. 
 
 Consider the \textit{WRLS}  problem \eqref{ls:eq12} and 
let  ${\bf r}=W(f-K^Tz).$ Thus, the minimization problem \eqref{ls:eq12} can be reformulated as the following \textit{SPP}:
\begin{align}\label{eq62}
	\widehat{\A}\bmatrix{{\bf r} \\ z}\triangleq \bmatrix{W^{-1} & K^{T}\\ K & -\lambda I_{ {m}}}\bmatrix{{\bf r} \\ z}=\bmatrix{f\\ {\bf 0}},
\end{align}
where $K$ is a  \textit{Toeplitz} or \textit{symmetric}-\textit{Toeplitz} matrix.

Since the weighting matrix $W$ and the regularization matrix $-\lambda I_{ {m}}$ 
 are not allowed to be perturbed, we consider  $(1,1)$-block  and $(2,2)$-block   has no perturbation. Let $[\widetilde{\bf r}^T,\, \widetilde{z}^T]^T$ be an approximate solution of of the system \eqref{eq62}, i.e.,  $\widetilde{z}$ be an approximate solution of the \textit{WRLS} problem. Then, we define structured \textit{BE} for the \textit{WRLS} problem as follows:
 \begin{align}
     \bm{\zeta}(\widetilde{z}):=\min_{\left(\D K, \D f\right)\in \mathcal{S}^{ls}}\left\|\bmatrix{[w_2\|\D K\|_F, & w_4\|\D f\|_2}\right\|_2,
 \end{align}
  where $$\mathcal{S}^{ls}:=\Bigg\{\left(\D K, \D f\right)\, : \bmatrix{W^{-1} & ( K+\D K)^T\\ K+\D K & -\lambda I_{ {m}}}\bmatrix{\widetilde{\r}\\ \widetilde{z}}=\bmatrix{f+\D f\\{\bf 0}}, \, \D K\in \{ {\T_{m\times n}}, \ST_n\}\Bigg\}.$$
 Before proceeding, we define $ \mathcal{X}_{ls}\in \R^{ {(n+m)}\times ( {2n+m-1})}$ and $ \mathcal{Z}_{ls}\in \R^{2n\times 2n}$ as follows:
\begin{align*}
\mathcal{X}_{ls}=\bmatrix{\frac{1}{w_2}\mathcal{G}_{\widetilde{z}}\mathfrak{D}_{t(K)}\mathfrak{D}_{{\bf t}_{mn}}^{-1}  &-\frac{1}{w_4}I_n \\
\frac{1}{w_2}\mathcal{K}_{\widetilde{\r}}\mathfrak{D}_{t(K)}\mathfrak{D}_{{\bf t}_{mn}}^{-1} & {\bf 0}}\quad \text{and}\quad\mathcal{Z}_{ls}=\bmatrix{\frac{1}{w_2}\mathcal{I}_{\widetilde{z}}\mathfrak{D}_{t(K)}\mathfrak{D}^{-1}_{\bf s}  &-\frac{1}{w_4}I_n \\
\frac{1}{w_2}\mathcal{I}_{\widetilde{\r}}\mathfrak{D}_{t(K)}\mathfrak{D}^{-1}_{\bf s} & {\bf 0}},
\end{align*}
 where $t(K)=\vec_{\T}(\Theta_K),$ $\widetilde{r}_f=f-W^{-1}\widetilde{\r}-K^T\widetilde{z},$ and $\widetilde{r}_g=\lambda \widetilde{z}-K\widetilde{\r}.$ 
  \begin{theorem}\label{th61}
  Let $\widetilde{z}$ be an approximate solution of the \textit{WRLS} problem \eqref{ls:eq12} with $K\in \{ {\T_{m\times n}}, \ST_n\}$.  {Let $\widetilde{r}_d=\bmatrix{\widetilde{r}_f^T,& \widetilde{r}_g^T}^T,$} then 
  \begin{enumerate}
      \item when $K\in  {\T_{m\times n}}$ and $\rank(\mathcal{X}_{ls})=\rank\left(\left[\mathcal{X}_{ls}\, \,\widetilde{r}_d \right]\right),$ we have
      \begin{align}\label{eq64}
          \bm{\zeta}(\widetilde{z})=\left\|\mathcal{X}_{ls}^{\dagger}{ \widetilde{r}_d}\right\|_2,
\end{align}
\item   when $K\in \ST_n$ and $\rank(\mathcal{Z}_{ls})=\rank\left(\left[\mathcal{Z}_{ls}\, \,\widetilde{r}_d \right]\right),$ we have \begin{align}\label{eq65}
          \bm{\zeta}(\widetilde{z})=\left\|\mathcal{Z}_{ls}^{\dagger}{ \widetilde{r}_d}\right\|_2.
\end{align}
  \end{enumerate}
 \end{theorem}
 \proof
 First, we consider $K\in  {\T_{m\times n}}.$
 Since, $W$ and $-\lambda I_{ {m}}$ are  not required to be perturbed, we take $w_1=0$ and $w_3=0.$ 
 Following the proof method of Theorem \ref{th41}, we obtain that $\left(\D K, \D f\right)\in \mathcal{S}^{ls},$ $\D K\in  {\T_{m\times n}}$ if and only if 
 \begin{eqnarray}\label{eq66}
     \mathcal{X}_{ls}\D \E_{ls}={ \widetilde{r}_d},
 \end{eqnarray} where $\D \E_{ls}=\bmatrix{w_2 {\mathfrak{D}_{{\bf t}_{mn}}}\vec_{\T}(\D K\,\odot\,\Theta_{K})\\ w_4\D f}.$ Thus, when $\rank(\mathcal{X}_{ls})=\rank\left(\left[\mathcal{X}_{ls}\, \,{ \widetilde{r}_d} \right]\right),$ the minimum norm solution of \eqref{eq66} is $\D \E^{ls}_{\min}=\mathcal{X}_{ls}^{\dagger}{ \widetilde{r}_d}.$  Hence,  the structured \textit{BE} ${\bm\zeta}(\widetilde{z})$ in \eqref{eq64} is attained. Similarly, for $K\in \ST_n\, (n=m),$ we can derive the structured \textit{BE} given in \eqref{eq65}. $\blacksquare$

Note that we can similarly obtain the structured \textit{BE} for the \textit{WRLS} when $K$ is circulant.

\section{Numerical experiments}\label{SEC:Numerical}
 In this section, we conduct a few numerical experiments to validate our findings. 
 All numerical experiments are conducted on MATLAB R2023b on an Intel(R) Core(TM) $i7$-$10700$ $CPU$, $ 2.90GHz,$ $ 16$ $GB$ memory with machine precision $\mu=2.2\times 10^{-16}.$   

{ \begin{exam}\label{EX:revision1}
   Consider a \textit{Toeplitz} structured \textit{SPP} where the coefficient matrices are given by:
    \begin{align*}
        &A=\bmatrix{10^{-6}& 0& 10^3 & 0\\ 10^8 &10^{-6} &0 & 10^3\\ 10 &10^8& 10^{-6}&0\\ 0 & 10 &10^8& 10^{-6} },\, ~B=\bmatrix{10^{-5}& 10^7& 0 & 0\\ 10^5 &10^{-5} &10^7 &0\\ 0 &10^5& 10^{-5}&10^7\\ 0 & 0 &10^5& 10^{-5}},\,\\
        &C=\bmatrix{0& 10^8 & -60& 0\\-0.5 &0& 10^8 &-60 \\ 0& -0.5 &0& 10^8 \\ 0& 0& -0.5 &0},\, f=\bmatrix{10^8\\ 0\\10^3\\0} ~\text{and}~g=\bmatrix{10^{-8}\\ 0\\0\\0}. 
    \end{align*}
  The approximate solution of the \textit{SPP}, computed using Gaussian elimination with partial pivoting (GEP) is  $ {\widetilde{u}}=[\x^{T},\, \y^{T}]^{T},$ where
$$\widetilde{x}=\bmatrix{6.0278\times 10^{3}\\
-1.0000\times 10^{4}\\
9.8995\times 10^{-3}\\
-9.9000\times 10^{7}} ~\text{and}~ \widetilde{y}=\bmatrix{-5.0378\times 10^{4}\\
1.0000\times 10^{3}\\
-8.8995\times 10^{-2}\\
9.9000\times 10^{6}}.$$

\noindent We choose $w_1=1/\|A\|_F,$ $W_2=1/\|B\|_F,$ $w_3=1/\|C\|_F,$ $w_4=1/1/\|f\|_2,$ and $w_5=1/\|g\|_2.$ Using  the formula in \eqref{eq14}, the unstructured \textit{BE} is $\eta( { {\widetilde{u}}})= 6.2617\times 10^{-18}.$ However, the obtained structured \textit{BEs} using Theorem \ref{th51} and Corollary \ref{coro51} are $\eta^{\mathcal{S}_2}(\x,\y)=2.1761\times 10^{-9}$ and $\eta^{\mathcal{S}_2}_{\tt sps}(\x,\y)=4.3070\times 10^{-5}$, respectively. We can observe that the $\eta( { {\widetilde{u}}})$ in the order of $\mathcal{O}(10^{-18}),$ whereas the structured \textit{BEs} are significantly larger. This demonstrates that the GEP for solving this tested \textit{Toeplitz} structured \textit{SPP} is backward stable but not strongly backward stable. That is, the computed approximate solution does not satisfy a nearby (sparsity preserving) \textit{Toeplitz} structured \textit{SPP}.
\end{exam}}
\begin{exam}\label{exam3}
Consider the \textit{Toeplitz} structured \textit{SPP} \eqref{eq11} with block matrices
\begin{align*}
    A=toeplitz({\bm a}_1, {\bm a}_2)\in \T_{ n\times n }, \quad  B=toeplitz({\bm b}_1, {\bm b}_2)\in \T_{ n\times n }, \quad  C=toeplitz({\bm c}_1, {\bm c}_2)\in \T_{ n\times n }, 
\end{align*}
 where ${\bm a}_1=sprand(n,1,0.4),$ ${\bm a}_2=randn(n,1),$ ${\bm b}_1=sprand(n,1,0.1),$ ${\bm b}_2=randn(n,1),$ ${\bm c}_1=sprand(n,1,0.1)$ and ${\bm c}_2=randn(n,1)$ so that ${\bm a}_1(1)={\bm a}_2(1), {\bm b}_1(1)={\bm b}_2(1)$ and ${\bm c}_1(1)={\bm c}_2(1).$  Moreover, we choose $f=randn(n,1)$ and $g=randn(n,1).$ Here, the notation $`randn(n,1)$' stands for normally distributed random vectors, and $`sprand(n,1,\omega)$' stands for   uniformly distributed sparse random vectors with density $\omega$ of size $n$. The symbol $toeplitz(\bm{a}_1,\bm{a}_2)$ denotes \textit{Toeplitz} matrix with $\bm{a}_1$ as its first column and $\bm{a}_2$ as its first row. We choose the parameters $w_i=1,$ for $i=1,2,\ldots,5.$
 \begin{figure}[]
					\centering
					\includegraphics[width=0.45\textwidth]{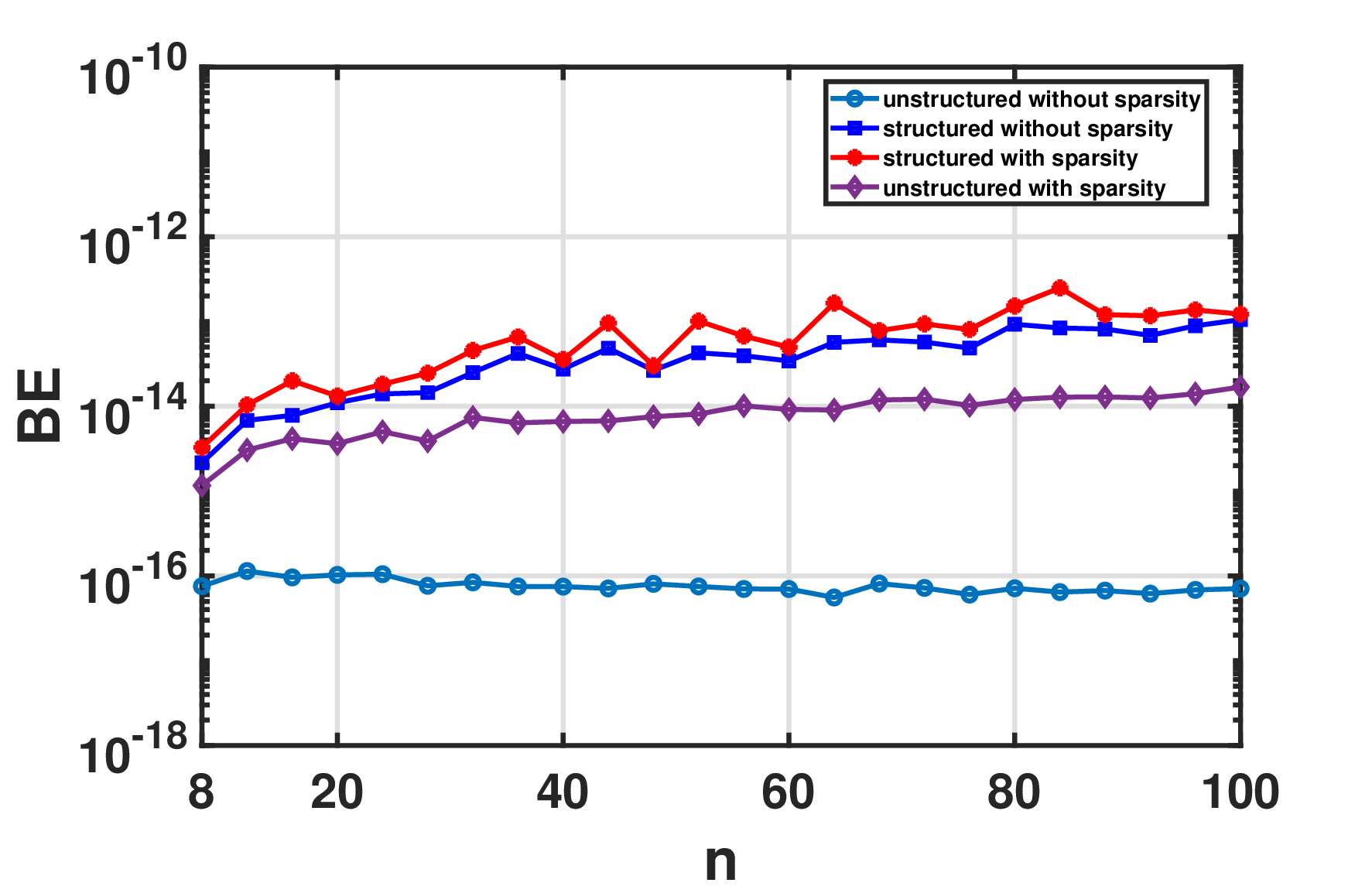}
		 \caption{Different structured and unstructured \textit{BEs} for $n=8:4:100$ .}
		\label{comparisonTeop}
	\end{figure}
 
We apply the \textit{GMRES} method \cite{gmres}  with the initial guess vector zero and tolerance $10^{-7}.$ Let $ {\widetilde{u}}=[\x^{T},\, \y^{T}]^{T}$ be the computed solution of the  \textit{SPP}. For $n=8:4:100$ in Figure \ref{comparisonTeop}, we plot the unstructured \textit{BE} $\eta( { {\widetilde{u}}})$ using the formula \eqref{eq14} (denoted as `\textbf{unstructured without sparsity}'), structured \textit{BE} $\eta_{\tt sps}^{\mathcal{S}_2}(\x, \y)$ (denoted as `\textbf{structured with sparsity}') using Theorem \ref{th41}, $\eta^{\mathcal{S}_2}(\x, \y)$ (denoted as `\textbf{structured without sparsity}') using Corollary \ref{coro41} and  $\eta_{\tt sps}(\x, \y)$ (denoted as `\textbf{unstructured without sparsity}') using Theorem \ref{th52}. From Figure \ref{comparisonTeop}, it can be observed that, for all values of $n,$ the unstructured \textit{BE} $\eta( {\widetilde{u}})$ around of order $\mathcal{O}(10^{-16})$  and all other \textit{BEs}  $\eta_{\tt sps}(\x, \y),\, \eta_{\tt sps}^{\mathcal{S}_2}(\x, \y), $ and $\eta^{\mathcal{S}_2}(\x, \y)$ are around of $\mathcal{O}(10^{-13})$ or less than that, which are very small.  {Therefore, the approximate solution computed using \textit{GMRES} method for each generated \textit{Toeplitz} structured \textit{SPP} effectively solves  a nearby perturbed unstructured linear system as well as a nearby perturbed (sparsity preserving) \textit{Toeplitz} structured \textit{SPP}.}
\end{exam}
  \begin{exam}\label{ex5}
    \begin{table}[ht!]
        \centering
        \caption{ Unstructured and structured \textit{BEs} for different values of $n$ for Example \ref{ex5}.}
        \resizebox{9cm}{!}{
       {  \begin{tabular}{ccccc}
        \toprule
         $n$ & $\eta( {\widetilde{u}})$&$\eta_{\tt sps}(\x, \y) $ &$\eta^{\mathcal{S}_3}(\x, \y) $& $\eta^{\mathcal{S}_3}_{\tt sps}(\x, \y)$ \\
             \midrule
           $8$& $1.5522e-16$	&$5.1676e-16$&	$1.3612e-15$&	$4.0978e-15 $\\
           \midrule
           $16$& $7.9177e-16$&	$2.6522e-15$&	$1.1505e-15$	&$8.3191e-15$\\
           \midrule
           $32$&$4.2992e-15$&	$1.1547e-13$	&$1.3009e-13$&	$6.2542e-13$\\
           \midrule 
           $64$ &$8.4321e-16$&	$4.8161e-14$&	$1.2379e-13$	&$4.2252e-1$3\\
           \midrule
           $128$ & $1.7815e-15$&	$3.4069e-14$&	$1.5579e-13$&	$6.1305e-13$\\
           \bottomrule
        \end{tabular}
        }}
        \label{tab:my_label}
    \end{table}
 { In this example, we consider the \textit{symmetric}-\textit{Toeplitz} structured \textit{SPP} \eqref{eq11} with the block matrices $A=I_n,$  $B=[b_{ij}]\in \ST_n,$  where $b_{ij}=\frac{1}{\sqrt{2\pi}}e^{-\frac{(i-j)^2}{2}},$ $i,j=1,\ldots, n,$ $C=-\mu I_n,$ $f=randn(n,1)\in \R^n,$ and $g={\bf 0}\in \R^n.$ We choose $\mu=0.01.$ To solve the \textit{SPP}, we use the  \textit{CNAGSOR} preconditioned \textit{GMRES} method \cite{circulantbased2021}. 
We choose $w_i,$ $i=1,\ldots, 4,$ as in Example \ref{EX:revision1} and $w_5=0.$ For the computed approximate solution  $ {\widetilde{u}}=[\x^T,\y^T]^T,$ we compute the unstructured \textit{BE} $\eta( { {\widetilde{u}}})$ using the formula \eqref{eq14}, structured \textit{BE} $\eta_{\tt sps}^{\mathcal{S}_3}(\x, \y)$ using Theorem \ref{th51}, and the structured \textit{BE} without preserving sparsity $\eta^{\mathcal{S}_3}(\x, \y)$ using Corollary \ref{coro51}. The computed values are reported in Table \ref{tab:my_label}. }

 {We observe that the structured \textit{BEs} $\eta_{\tt sps}^{\mathcal{S}_3}(\x, \y)$ and $\eta^{\mathcal{S}_3}(\x, \y)$ are all most all cases remains one or two order larger than the unstructured ones and remains within an order of $\mathcal{O}(10^{-13}).$ Hence, we can conclude that the approximate solution \( {\widetilde{u}}\) obtained using the \textit{CNAGSOR}-preconditioned \textit{GMRES} method for the tested \textit{SPPs} serves as an exact solution to a nearly perturbed \textit{symmetric}-\textit{Toeplitz} structured \textit{SPP}, while preserving the sparsity pattern of the original problem.} 
  \end{exam}  
 \section{Conclusion}\label{SEC:Conclusion}
 In this paper, we investigate the structured \textit{BEs} for \textit{ {\textit{circulant}}}, \textit{Toeplitz}, and \textit{symmetric}-\textit{Toeplitz} structured \textit{SPPs} with and without preserving the sparsity pattern of block matrices. Additionally, we provide minimal perturbation matrices for which an approximate solution becomes the exact solution of a nearly perturbed \textit{SPP}, which preserves the inherent block structure and sparsity pattern of the original \textit{SPP}. Furthermore, unstructured \textit{BE} is obtained when the block matrices of the \textit{SPPs} only preserve the sparsity pattern. Our obtained results are used to derive structured \textit{BE} for \textit{WRLS} problems with \textit{Toeplitz} or \textit{symmetric}-\textit{Toeplitz} coefficient matrices. Numerical experiments are performed to validate our theoretical findings and to examine the backward stability and the strong backward stability of numerical algorithms to solve structured \textit{SPPs}. 
 \section*{Acknowledgments}
Pinki Khatun acknowledges the Council of Scientific $\&$ Industrial Research (CSIR) in New Delhi, India, for the financial support provided to her in the form of a fellowship (File no. 09/1022(0098)/2020-EMR-I). 
 {The authors express their gratitude to the anonymous reviewers for their insightful suggestions, which have significantly enhanced the quality of this paper.}

 \bibliography{BEreference}

\begin{thebibliography}{28}
\providecommand{\natexlab}[1]{#1}
\providecommand{\url}[1]{\texttt{#1}}
\expandafter\ifx\csname urlstyle\endcsname\relax
  \providecommand{\doi}[1]{doi: #1}\else
  \providecommand{\doi}{doi: \begingroup \urlstyle{rm}\Url}\fi

\bibitem[Ahmad and Kanhya(2020)]{prince2020}
S.~S. Ahmad and P.~Kanhya.
\newblock Structured perturbation analysis of sparse matrix pencils with
  {$s$}-specified eigenpairs.
\newblock \emph{Linear Algebra Appl.}, 602:\penalty0 93--119, 2020.

\bibitem[Ahmad and Kanhya(2021)]{prince2021}
S.~S. Ahmad and P.~Kanhya.
\newblock Backward error analysis and inverse eigenvalue problems for {H}ankel
  and {S}ymmetric-{T}oeplitz structures.
\newblock \emph{Appl. Math. Comput.}, 406:\penalty0 126288, 15, 2021.

\bibitem[Balani and Hajarian(2023)]{weightedLS2023}
F.~B. Balani and M.~Hajarian.
\newblock A new block preconditioner for weighted {T}oeplitz regularized
  least-squares problems.
\newblock \emph{Int. J. Comput. Math.}, 100\penalty0 (12):\penalty0 2241--2250,
  2023.

\bibitem[Benzi et~al.(2005)Benzi, Golub, and Liesen]{Benzi2005}
M.~Benzi, G.~H. Golub, and J.~Liesen.
\newblock Numerical solution of saddle point problems.
\newblock \emph{Acta Numer.}, 14:\penalty0 1--137, 2005.

\bibitem[Bolten et~al.(2023)Bolten, Donatelli, Ferrari, and
  Furci]{SymbolCirculant}
M.~Bolten, M.~Donatelli, P.~Ferrari, and I.~Furci.
\newblock Symbol based convergence analysis in multigrid methods for saddle
  point problems.
\newblock \emph{Linear Algebra Appl.}, 671:\penalty0 67--108, 2023.

\bibitem[Bunch(1987)]{strongweak}
J.~R. Bunch.
\newblock The weak and strong stability of algorithms in numerical linear
  algebra.
\newblock \emph{Linear Algebra Appl.}, 88/89:\penalty0 49--66, 1987.

\bibitem[Bunch et~al.(1989)Bunch, Demmel, and Van~Loan]{strongstab}
J.~R. Bunch, J.~W. Demmel, and C.~F. Van~Loan.
\newblock The strong stability of algorithms for solving symmetric linear
  systems.
\newblock \emph{SIAM J. Matrix Anal. Appl.}, 10\penalty0 (4):\penalty0
  494--499, 1989.

\bibitem[Chang et~al.(2009)Chang, Paige, and Titley-Peloquin]{stopping2009}
X.-W. Chang, C.~C. Paige, and D.~Titley-Peloquin.
\newblock Stopping criteria for the iterative solution of linear least squares
  problems.
\newblock \emph{SIAM J. Matrix Anal. Appl.}, 31\penalty0 (2):\penalty0
  831--852, 2009.

\bibitem[Chen et~al.(2012)Chen, Li, Chen, and Liu]{be2012LAA}
X.~S. Chen, W.~Li, X.~Chen, and J.~Liu.
\newblock Structured backward errors for generalized saddle point systems.
\newblock \emph{Linear Algebra Appl.}, 436\penalty0 (9):\penalty0 3109--3119,
  2012.

\bibitem[Elman et~al.(2005)Elman, Silvester, and Wathen]{Elman2005}
H.~C. Elman, D.~J. Silvester, and A.~J. Wathen.
\newblock \emph{Finite Elements and Fast Iterative Solvers: With Applications
  in Incompressible Fluid Dynamics}.
\newblock Oxford University Press, New York, 2005.

\bibitem[Fessler and Booth(1999)]{imagereconstruction}
J.~A. Fessler and S.~D. Booth.
\newblock Conjugate-gradient preconditioning methods for shift-variant pet
  image reconstruction.
\newblock \emph{IEEE Transactions on Image Processing}, 8\penalty0
  (5):\penalty0 688--699, 1999.

\bibitem[Greif et~al.(2014)Greif, Moulding, and Orban]{OPTM2014}
C.~Greif, E.~Moulding, and D.~Orban.
\newblock Bounds on eigenvalues of matrices arising from interior-point
  methods.
\newblock \emph{SIAM J. Optim.}, 24\penalty0 (1):\penalty0 49--83, 2014.

\bibitem[Gulliksson et~al.(2002)Gulliksson, Jin, and Wei]{LSproblem2002}
M.~Gulliksson, X.-Q. Jin, and Y.-M. Wei.
\newblock Perturbation bounds for constrained and weighted least squares
  problems.
\newblock \emph{Linear Algebra Appl.}, 349:\penalty0 221--232, 2002.

\bibitem[Higham(2002)]{higham2002}
N.~J. Higham.
\newblock \emph{Accuracy and Stability of Numerical Algorithms}.
\newblock SIAM, Philadelphia, PA, second edition, 2002.

\bibitem[Jain(1989)]{image1989}
A.~K. Jain.
\newblock \emph{Fundamentals of digital image processing}.
\newblock Prentice-Hall, Inc., 1989.

\bibitem[Li and Liu(2004)]{BEKKT2004}
X.~Li and X.~Liu.
\newblock Structured backward errors for structured {KKT} systems.
\newblock \emph{J. Comput. Math.}, 22\penalty0 (4):\penalty0 605--610, 2004.

\bibitem[Ma(2017)]{be2017ma}
W.~Ma.
\newblock On normwise structured backward errors for the generalized saddle
  point systems.
\newblock \emph{Calcolo}, 54\penalty0 (2):\penalty0 503--514, 2017.

\bibitem[Meng et~al.(2022)Meng, He, and Miao]{be2022lma}
L.~Meng, Y.~He, and S.-X. Miao.
\newblock Structured backward errors for two kinds of generalized saddle point
  systems.
\newblock \emph{Linear Multilinear Algebra}, 70\penalty0 (7):\penalty0
  1345--1355, 2022.

\bibitem[Rigal and Gaches(1967)]{Rigal1967}
J.~L. Rigal and J.~Gaches.
\newblock On the compatibility of a given solution with the data of a linear
  system.
\newblock \emph{J. Assoc. Comput. Mach.}, 14:\penalty0 543--548, 1967.

\bibitem[Saad and Schultz(1986)]{gmres}
Y.~Saad and M.~H. Schultz.
\newblock G{MRES}: a generalized minimal residual algorithm for solving
  nonsymmetric linear systems.
\newblock \emph{SIAM J. Sci. Statist. Comput.}, 7\penalty0 (3):\penalty0
  856--869, 1986.

\bibitem[Sun(1999)]{Sun1999}
J.~G. Sun.
\newblock Structured backward errors for {KKT} systems.
\newblock \emph{Linear Algebra Appl.}, 288\penalty0 (1-3):\penalty0 75--88,
  1999.

\bibitem[Wang et~al.(2018)Wang, Wei, and Qiao]{GINVERSE}
G.~Wang, Y.~Wei, and S.~Qiao.
\newblock \emph{Generalized Inverses: Theory and Computations}.
\newblock Springer, Singapore, second edition, 2018.

\bibitem[Wright(1992)]{OPTM1992}
M.~H. Wright.
\newblock Interior methods for constrained optimization.
\newblock \emph{Acta Numer.}, 1:\penalty0 341--407, 1992.

\bibitem[Xiang and Wei(2007)]{be2007wei}
H.~Xiang and Y.~Wei.
\newblock On normwise structured backward errors for saddle point sytems.
\newblock \emph{SIAM J. Matrix Anal. Appl.}, 29\penalty0 (3):\penalty0
  838--849, 2007.

\bibitem[Zeng(2021)]{circulantbased2021}
M.~L. Zeng.
\newblock A circulant-matrix-based new accelerated {GSOR} preconditioned method
  for block two-by-two linear systems from image restoration problems.
\newblock \emph{Appl. Numer. Math.}, 164:\penalty0 245--257, 2021.

\bibitem[Zhang and Su(2012)]{PEVP2012}
K.~Zhang and Y.~Su.
\newblock Structured backward error analysis for sparse polynomial eigenvalue
  problems.
\newblock \emph{Appl. Math. Comput.}, 219\penalty0 (6):\penalty0 3073--3082,
  2012.

\bibitem[Zheng and Lv(2020)]{BE2020BING}
B.~Zheng and P.~Lv.
\newblock Structured backward error analysis for generalized saddle point
  problems.
\newblock \emph{Adv. Comput. Math.}, 46\penalty0 (2):\penalty0 34--27, 2020.

\bibitem[Zhu et~al.(2018)Zhu, Qi, and Zhang]{SYMTOEP2017}
M.~Z. Zhu, Y.~E. Qi, and G.~F. Zhang.
\newblock On local circulant and residue splitting iterative method for
  {T}oeplitz-structured saddle point problems.
\newblock \emph{IAENG Int. J. Appl. Math.}, 48\penalty0 (2):\penalty0 221--227,
  2018.

\end{thebibliography}
\bibliographystyle{abbrvnat}
 \end{document}